\documentclass{amsart}
\usepackage{amssymb, latexsym}
\newcommand{\half}{\frac{1}{2}}
\newcommand{\third}{\frac{1}{3}}

\newcommand{\e}{\epsilon}

\newcommand{\lam}{\lambda}
\newcommand{\Lam}{\Lambda}

\newcommand{\R}{{\mathbb R}}

\newcommand{\C}{{\mathbb C}}
\newcommand{\T}{{\mathbb T}}

\newcommand{\Z}{{\mathbb Z}}
\newcommand{\N}{{\mathbb N}}

\newcommand{\wt}{\widetilde}

\newcommand{\uhat}{\widehat{u}}

\newcommand{\Lp}[2]{{\left\| #2 \right\|}_{L^ #1 }}

\newcommand{\Hsup}[2]{{{\left\| #2 \right\|}_{{H^#1}}}}

\newcommand{\X}[3]{{{\left\| #3 \right\|}_{{X_{#1,#2}}}}}
\newcommand{\Xd}[3]{{{\left\| #3 \right\|}_{{X^\delta_{#1,#2}}}}}

\newcommand{\Xsb}{{X_{s,b}}}

\newcommand{\xionetofour}{\xi_1 \xi_2 \xi_3 \xi_4 }
\newcommand{\bint}{\int\limits_{\xi_1 + \xi_2 = 0}}

\newtheorem{theorem}{Theorem}
\theoremstyle{definition}
\newtheorem{definition}{Definition}
\theoremstyle{remark}
\newtheorem{remark}{Remark}
\theoremstyle{proposition}
\newtheorem{proposition}{Proposition}
\theoremstyle{lemma}
\newtheorem{lemma}{Lemma}
\theoremstyle{corollary}

\numberwithin{equation}{section}
\numberwithin{lemma}{section}
\numberwithin{remark}{section}

\begin{document}

\title{Sharp Global well-posedness for KdV and modified KdV on $\R$ and $\T$}

%%%%%%%%%%%%%%%%%%%%%%%%%%%%%%%%%%%%%%%%%%%%%%%%%%%%%%%%%%%%%%%%%
%                                                               %
%                          The I-Team                           %
%                                                               %
%%%%%%%%%%%%%%%%%%%%%%%%%%%%%%%%%%%%%%%%%%%%%%%%%%%%%%%%%%%%%%%%%

\author{J.~Colliander}
\thanks{J.C. is supported in part by N.S.F. Grant DMS 0100595.}
\address{University of Toronto}
\author{M.~Keel}
\thanks{M.K. is supported in part by N.S.F. Grant DMS 
                         9801558}
\address{University of Minnesota}
\author{G.~Staffilani}
\thanks{G.S. is supported in part by N.S.F. Grant DMS 9800879 and by a grant
from the Sloan Foundation.}
\address{Brown University}
\author{H.~Takaoka}
\address{Hokkaido University}
\thanks{H.T. is supported in part by J.S.P.S. Grant No. 13740087.}
\author{T.~Tao}
\thanks{T.T. is a Clay Prize Fellow and is supported in part by grants
from the Packard and Sloan Foundations.}
\address{University of California, Los Angeles}

\subjclass{35Q53, 42B35, 37K10}
\keywords{Korteweg-de Vries equation, nonlinear dispersive equations,
bilinear estimates, 
multilinear harmonic analysis} 

\begin{abstract}
The initial value problems for the Korteweg-de Vries (KdV) and
modified KdV (mKdV) equations under periodic and decaying boundary conditions
are considered. These initial value problems are shown to be globally
well-posed in all $L^2$-based Sobolev spaces $H^s$ 
where local well-posedness is presently known, apart from the 
$H^{\frac{1}{4}} (\R )$ endpoint for mKdV. The
result for KdV relies on a new method for constructing almost conserved
quantities using multilinear harmonic analysis and the available local-in-time
theory. Miura's transformation is used to show that global well-posedness
of modified KdV is implied by global well-posedness of the standard
KdV equation.

\end{abstract}

\date{3 October 2001}

\maketitle

\tableofcontents

\section{Introduction}

The initial value problem for the Korteweg-de Vries (KdV) equation, 
\begin{equation}
  \label{KdV}
  \left\{
   \begin{matrix}
    \partial_t u + \partial_x^3 u + \half \partial_x u^2 =0,& u: \R \times 
[0,T] \longmapsto \R ,
         \\
     u( 0) = \phi \in H^s ( \R ),
   \end{matrix}
\right.
\end{equation}
has been shown to be locally well-posed (LWP) for $s> - \frac{3}{4}.$
Kenig, Ponce and Vega \cite{KPVBilin} extended the local-in-time
analysis of Bourgain \cite{B1}, valid for $s \geq 0$, to the range 
$s > - \frac{3}{4}$ by constructing the solution of \eqref{KdV} on a 
time interval $[0, \delta]$ with $\delta$ depending upon 
${{\| \phi \|}_{H^s (\R )}}$. 
Earlier results can be found in \cite{BonaSmith75}, 
\cite{Kato81}, \cite{GinibreTsutsumiVelo}, \cite{KPVCPAM}, \cite{Cohen79}.  
We prove here that these solutions exist for $t$ in an arbitrary time interval
$[0,T]$ thereby establishing global well-posedness (GWP) of \eqref{KdV} in the
full range $s > - \frac{3}{4}.$
The corresponding periodic $\R$-valued initial value problem for KdV
\begin{equation}
  \label{TKdV}
  \left\{
   \begin{matrix}
    \partial_t u + \partial_x^3 u + \half \partial_x u^2 =0,& u: \T \times 
[0,T] \longmapsto \R
         \\
     u( 0) = \phi \in H^s ( \T )
   \end{matrix}
\right.
\end{equation}
is known \cite{KPVBilin} to be locally well-posed for $s \geq - \half$.
These local-in-time solutions are also shown to exist on an arbitrary
time interval. Bourgain established \cite{BKdVMeasures} global well-posedness
of \eqref{TKdV} for initial data having (small) bounded Fourier transform.
The argument in \cite{BKdVMeasures} uses the complete integrability of KdV.
Analogous globalizations of the best known local-in-time theory for the 
focussing and defocussing modified KdV (mKdV) equations ($u^2$ in 
\eqref{KdV}, \eqref{TKdV} replaced by $-u^3$ and $u^3$, respectively) are
also obtained in the periodic $(s \geq \half)$ and real line ($s \geq 
\frac{1}{4})$ settings.

It appears likely that the local-in-time theory globalized here is sharp,
at least in the scale of $L^2$-based Sobolev spaces $H^s$. Indeed,
recent examples \cite{KPVCounter} of Kenig, Ponce and Vega reveal that 
focussing mKdV is ill-posed for $s < \frac{1}{4}$ and that $\C$-valued
KdV ($u: \R \times [0,T] \longmapsto \C$) is illposed for $s < -\frac{3}{4}$.
(The local theory in \cite{KPVBilin} adapts easily to the $\C$-valued
situation.) The fundamental bilinear estimate used to prove the local
well-posedness result on the line was shown to fail for $s \leq -\frac{3}{4}$
by Nakanishi, Takaoka and Tsutsumi \cite{NTT}. A similar failure
for $s< - \half$ occurs \cite{KPVBilin} in the periodic problem.
However, a definitive understanding of the $\R$-valued initial
value problem for KdV or defocussing mKdV in the $H^s$ scale has not yet
been obtained. Furthermore, no ill-posedness result for any defocussing
dispersive problem is known.  

\subsection{GWP below the conservation law}

$\R$-valued solutions of KdV satisfy $L^2$ conservation: ${{\| u(t) \|}_{L^2}}
= {{\| \phi \|}_{L^2}}$. Consequently, a local well-posedness result
with the existence lifetime determined by the size of the initial data
in $L^2$ may be iterated to prove global well-posedness of KdV for $L^2$
data \cite{B1}. What happens to solutions of KdV which evolve from initial
data which are less regular than $L^2$? Bourgain observed, in a context 
\cite{B4} concerning very smooth solutions, that the nonlinear Duhamel term
may be more smooth than the initial data. This observation was exploited 
\cite{B4}, using a decomposition of the evolution of the high and 
low frequency parts of the initial data, to prove polynomial-in-time 
bounds for global solutions of certain nonlinear Schr\"odinger (NLS) and 
nonlinear wave (NLW) equations. In \cite{BRefine}, Bourgain introduced a 
general high/low frequency decomposition argument to prove 
that certain NLS and NLW equations were globally well-posed below 
$H^1$, the natural regularity associated with the conserved Hamiltonian.
Subsequently, Bourgain's high/low trick has been applied to prove
global well-posedness below the natural regularity of the conserved 
quantity in various settings \cite{FLP}, \cite{KeeTao98},  \cite{KeelTaoMKG},
\cite{TzvGWPKPII99}, \cite{TakGWPKPII99}, \cite{KPVWaveGWP}, including
KdV \cite{CST99} on the line.

We summarize the adaptation \cite{CST99} 
of the high/low trick to treat \eqref{KdV} 
below $L^2$. The task is to construct the solution of \eqref{KdV} evolving from
initial data $\phi \in H^s ( \R )$ for $s_0 < s < 0$
with $-3/4 \ll  s_0 \lesssim 0$. Split the data $\phi = \phi_0 + \psi_0$
with ${\widehat{\phi_0}} (k) =  \chi_{[-N,N]} (k) {\widehat{\phi}} (k) $,
where $N = N(T)$ is a parameter to be determined. The low frequency part
$\phi_0$ of $\phi$ is in $L^2 (\R )$ (in fact $\phi_0 \in H^s$ for all $s$) 
with a big norm
while the high frequency part $\psi_0$ is the tail of an $H^s (\R )$ 
function and
is therefore small (with large $N$) in $H^\sigma (\R)$ for any $\sigma < s$.
The low frequencies are evolved according to KdV: $\phi_0 \longmapsto 
u_0 (t).$ The high frequencies
evolve according to a ``difference equation'' which is selected so that
the sum of the resulting high frequency evolution, $\psi_0 \longmapsto 
v_0 (t)$ and the low frequency evolution solves \eqref{KdV}. The key
step is to decompose $v_0 (t) = S(t) \psi_0 + w_0 (t)$, where $S(t)$
is the solution operator of the Airy equation, and to prove that $w_0 \in
L^2 ( \R)$ and has a small (depending upon $N$) $L^2$ norm. Then an iteration
of the local-in-time theory advances the solution to a long (depending on $N$)
time interval. An appropriate choice of $N$ completes the construction.

The nonlinear Duhamel term for the ``difference equation'' is 
\begin{equation*}
w_0 (t) = \int_0^t S(t - t' ) ( [v_0^2 (t') + 2 u_0 (t') v_0( t') ])dt' .
\end{equation*}
The local well-posedness machinery \cite{B1}, \cite{KPVBilin} allows
us to prove that $w_0 (t) \in L^2 (\R )$ if we prove the {\it{extra smoothing 
bilinear estimate}}
\begin{equation}
\label{extrasmoothing}
  \X 0 {b-1} {\partial_x (u v) } \lesssim
\X s b u  \X s b v , ~{s < 0, b= \half +},
\end{equation}
with the space $\Xsb $ defined below (see \eqref{Xsbdefined}). Indeed, this
estimate, combined with the known inhomogeneous estimate in $\Xsb$, 
shows that for $u$ and $v$ having $H^s$ spatial
regularity, the resulting nonlinear Duhamel term is in $X_{0,b}$, and hence
has $L^2$ spatial regularity at each $t$. The estimate \eqref{extrasmoothing}
fails for $s < - \frac{3}{8}$ and this places an intrinsic limitation on
how far the high/low frequency decomposition technique may be used to
extend GWP below $L^2$. For functions $u,~v$ such that 
$\widehat{u}, ~\widehat{v}$ are supported
outside $\{ |k | \leq 1 \}$, \eqref{extrasmoothing} is valid in the range
$-\frac{3}{8} < s$ \cite{CST99}, \cite{CKSTTKdVFirstGen}. 
Low frequencies in $u$ and $v$ degrade the strength of 
\eqref{extrasmoothing} and created certain technical difficulties in
carrying out the high/low strategy to treat KdV in \cite{CST99}. In 
principle, there should be no difficulty arising from low frequencies
in globalizing rough solutions of KdV since the difficulties associated
with low regularity should be due to too much ``mass'' at high frequencies.
We showed that the low frequency issue may indeed be circumvented in 
\cite{CKSTTKdVFirstGen} by proving \eqref{KdV} is GWP in $H^s ( \R ),~
s> - \frac{3}{10}$.

\subsection{The operator $I$ and almost conserved quantitites}

Global well-posedness follows from (an iteration of) local well-posedness 
(results) provided the successive local-in-time existence intervals cover
an arbitrary time interval $[0,T]$. The length of the local-in-time existence
interval is controlled from below by the size of the initial data in an
appropriate norm. A natural approach to global well-posedness in $H^s$
is to establish upper bounds on $\Hsup s {u(t)} $ for solutions $u(t)$
which are strong enough to prove $[0,T]$ may be covered by iterated
local existence intervals. We establish appropriate upper bounds to carry
out this general strategy by constructing {\it{almost conserved quantities}}
and rescaling. The rescaling exploits the {\it{subcritical}} nature of
the KdV initial value problem (but introduces technical issues in the
treatment of the periodic problem). The almost conserved quantities
are motivated by the following discussion of the $L^2$ conservation property
of solutions of KdV.

Consider the following Fourier 
proof{\footnote{This argument was known previously \cite{RauchPrivate}; similar
arguments appear in
\cite{JMR95}.}} that $\Lp 2 {u(t)} = \Lp 2 \phi ~\forall t \in \R$. By 
Plancherel,
\begin{equation*}
{{\| u(t) \|}_{L^2}^2} = \int \widehat{u} ( \xi )
\overline{\widehat{u}} ( \xi) d \xi,
\end{equation*}
where
\begin{equation*}
  \widehat{u} ( \xi ) = \int e^{- i x \xi } u(x) dx
\end{equation*}
is the (spatial) Fourier transform.
Fourier transform properties imply
\begin{equation*}
\int \widehat{u} (\xi ) \overline{\widehat{u}} ( \xi) d \xi =
\int \widehat{u} ( \xi) {\widehat{\overline{u}}} ( - \xi ) d \xi
= \int\limits_{\xi_1 + \xi_2 = 0}  \widehat{u} ( \xi_1 ) 
\widehat{\overline{u}} (\xi_2 ).
\end{equation*}
Since we are assuming $u$ is $\R$-valued, we may replace 
$\widehat{\overline{u}}( \xi_2) $ by $\widehat{u} ( \xi_2 )$. Hence,
\begin{equation*}
{{\| u(t) \|}_{L^2}^2} = \int\limits_{\xi_1 + \xi_2 = 0} 
\widehat{u} ( \xi_1 ) \widehat{u} (\xi_2 ).
\end{equation*}
We apply $\partial_t $, use symmetry, and the equation to find
\begin{equation*}
  \partial_t {{( {{\| u(t) \|}_{L^2}^2 })}} = 
2i \int\limits_{\xi_1 + \xi_2 =0 } \xi_1^3 \widehat{u} ( \xi_1 ) \widehat{u}
(\xi_2) - i \int\limits_{\xi_1 + \xi_2 =0 } \xi_1 \widehat{ u^2} (\xi_1 )
\widehat{u} ( \xi_2 ).
\end{equation*}
The first expression is symmetric under the interchange of $\xi_1$ and $\xi_2$
so $\xi_1^3$ may be replaced by $\half ( \xi_1^3 + \xi_2^3 )$. Since we
are integrating on the set where $\xi_1 + \xi_2 = 0$, the integrand is zero
and this term vanishes. Calculating $\widehat{u^2} ( \xi ) = \int\limits_{\xi 
= \xi_1 + \xi_2 } \widehat{u} ( \xi_1 ) \widehat{u} ( \xi_2 )$, the remaining
term may be rewritten
\begin{equation}
\label{nonlinLTwo}
  -i \int_{\xi_1 + \xi_2 + \xi_3 = 0} [ \xi_1 + \xi_2 ]
\uhat (\xi_1 ) \uhat ( \xi_2 ) \uhat ( \xi_3 ).
\end{equation}
On the set where $\xi_1 + \xi_2 + \xi_3 = 0 $, $\xi_1 + \xi_2 = - \xi_3 $
which we symmetrize to replace $\xi_1 + \xi_2$ in \eqref{nonlinLTwo} by
$-\frac{1}{3} ( \xi_1 + \xi_2 + \xi_3 )$ and this term vanishes as well.
Summarizing, we have found that $\R$-valued solutions $u(t)$ of KdV
satisfy
\begin{equation}
  \label{LTwoFourier}
  \partial_t {{( {{\| u(t) \|}_{L^2}^2} )}} = -i \bint (\xi_1^3 + \xi_2^3 )
\uhat (\xi_1 ) \uhat (\xi_2) + \frac{i}{3} \int\limits_{\xi_1 + \xi_2 
+ \xi_3 =0} (\xi_1 + \xi_2 + \xi_3 ) \uhat (\xi_1 ) \uhat (\xi_2) \uhat (\xi_3)
\end{equation}
and both integrands on the right-side vanish.

We introduce the (spatial) Fourier multiplier operator $Iu$ defined via
\begin{equation*}
  {\widehat{Iu}} (\xi ) = m(\xi ) \uhat (\xi )
\end{equation*}
with an arbitrary $\C$-valued multiplier $m$. A formal imitation of the
Fourier proof of $L^2$-mass conservation above reveals that for $\R$-valued
solutions of KdV we have
\begin{equation}
  \label{IulikeLTwo}
  \partial_t  {{( {{\| I u(t) \|}_{L^2}^2} )}} = - \frac{i}{2}
\bint [ m( \xi_1 ) \overline{m} ( \xi_1 ) + m( \xi_2 ) \overline{m} ( \xi_2 )]
\{ \xi_1^3 + \xi_2^3 \} ~\uhat ( \xi_1 ) \uhat ( \xi_2 )
\end{equation}
\begin{equation*}
  + \frac{i}{6} \int\limits_{\xi_1 + \xi_2 + \xi_3 = 0}
\sum\limits_{j=1}^3 [ m( - \xi_j ) \overline{m} ( - \xi_j ) + m( \xi_j )
\overline{m} ( \xi_j )] \xi_j ~\uhat (\xi_1 ) \uhat (\xi_2 ) \uhat ( \xi_3 ).
\end{equation*}
The term arising from the dispersion cancels since $ \xi_1^3 + \xi_2^3 =0$
on the set where $\xi_1 + \xi_2 =0$. The remaining trilinear term can
be analyzed under various assumptions on the multiplier $m$ giving insight
into the time behavior of ${{\| I u (t) \|}_{L^2}}$. Moreover, the 
flexibility in our choice of $m$ may allow us to observe how 
the conserved $L^2$
mass is moved around in frequency space during the KdV evolution.

\begin{remark} 
Our use of the multiplier $m$ to localize the $L^2$ mass in frequency
space is analagous to the use of cutoff functions
to spatially 
localize the conserved density on the spatial side. In that setting, 
the underlying conservation law $\partial_t ( ~{\mbox{conserved density}}~ ) 
+ \partial_x ( ~{\mbox{flux}}~ ) = 0$ is multiplied by a cutoff function.
The localized flux term is no longer a perfect derivative and is 
then estimated, sometimes under an appropriate choice of the cutoff, to
obtain bounds on the spatially localized energy. 
\end{remark}

Consider now the problem of proving well-posedness of \eqref{KdV} or
\eqref{TKdV}, with $s< 0$, on an arbitrary time interval $[0,T]$. We
define a spatial Fourier multiplier operator $I$ which acts like the
identity on low frequencies and like the $H^s$-norm on high frequencies
by choosing a smooth monotone multiplier satisfying
\begin{equation*}
m(\xi ) = \left\{ \begin{matrix}
1, & |\xi | < N \\
N^{-s} |\xi |^s , & |\xi | > 2N .
\end{matrix}
\right.
\end{equation*} 
The parameter $N$ marks the transition from low to high frequencies. 
When $N=1$, the operator $I$ is essentially the integration (since $s<0$)
operator $D^s$. When $N= \infty$, $I$ acts like the identity operator.
Note that ${{\| I \phi \|}_{L^2}}$ is bounded if $\phi \in H^s $. We prove a variant local well-posedness result which shows the length of the local existence interval $[0, \delta]$ for \eqref{KdV} or \eqref{TKdV} may be bounded 
from below by
${{\| I \phi \|}_{L^2}^{-\alpha}},
~\alpha > 0$ for an appropriate range of the parameter $s$. 
The basic idea is then to bound the trilinear term in \eqref{IulikeLTwo}
to prove, for a particular small $\beta > 0$, that 
\begin{equation}
\label{basicidea}
\sup_{t \in [0, \delta]} {{\| I u(t) \|}_{L^2}} \leq {{\| I u(0) \|}_{L^2}}
+ c N^{-\beta } {{\| I u (0 ) \|}_{L^2}^3}. 
\end{equation}

If $N$ is huge, \eqref{basicidea}
shows there is at most a tiny increment in ${{\| I u(t) \|}_{L^2}}$ as
$t$ evolves from $0$ to $\delta $. An iteration of the local theory under
appropriate parameter choices gives global well-posedness in $H^s$
for certain $s < 0$.

The strategy just described is enhanced with two extra ingredients:
a multilinear correction technique and rescaling. The correction technique
shows that, up to errors of smaller order in $N$, the trilinear term
in \eqref{IulikeLTwo} may be replaced by a quintilinear term improving
\eqref{basicidea} to
\begin{equation}
\label{betteridea}
\sup_{t  \in [0, \delta]} {{\| I u(t) \|}_{L^2}} \leq {{\| I u(0) \|}_{L^2}}
+ c N^{-3 - \frac{3}{4} + \epsilon } {{\| I u (0 ) \|}_{L^2}^5}, 
\end{equation}
where $\epsilon$ is tiny. The rescaling argument reduces matters to initial
data $\phi$ of fixed size: ${{\| I \phi \|}_{L^2}} \thicksim \epsilon_0 \ll 1$.
In the periodic setting, the rescaling we use forces us to track the
dependence upon the spatial period in the local well-posedness theory 
\cite{B1}, \cite{KPVBilin}.

The main results obtained here are:

\begin{theorem}
The initial value problem \eqref{KdV} is globally well-posed for 
initial data $\phi \in H^s (\R ),~ s> - \frac{3}{4}.$ 
\end{theorem}

\begin{theorem}
The periodic initial value problem \eqref{TKdV} is globally well-posed
for initial data $\phi \in H^s (\T ),~ s \geq -\half $.
\end{theorem}

\begin{theorem}
The initial value problem for modified KdV \eqref{mKdVivp} (focussing
or defocussing) is globally well-posed for initial data $\phi \in H^s
( \R ), ~ s > \frac{1}{4}$.
\end{theorem}

\begin{theorem}
The periodic initial value problem for modified KdV (focussing or 
defocussing) is globally well-posed for initial data $\phi \in H^s ( \T ), ~
s \geq \half $.
\end{theorem}

Our results here and elsewhere \cite{CKSTTDNLS1}, \cite{CKSTTDNLS2},
\cite{CKSTTGKdV} suggest that, in the absence of a mechanism for blow-up,
local well-posedness implies global well-posedness in subcritical dispersive 
initial value problems. In particular, we believe our methods will extend
to prove GWP of mKdV in $H^{\frac{1}{4}} ( \R )$.

The infinite dimensional symplectic nonsqueezing machinery developed by 
S. Kuksin \cite{KuksinSqueeze} identifies $H^{-\half} (\T) $ as the Hilbert 
Darboux (symplectic) phase space for KdV. We anticipate that Theorem
3 will be useful in adapting these ideas to
the KdV context. The main remaining issue is an approximation of the KdV flow
using finite dimensional Hamiltonian flows analagous to that obtained
by Bourgain \cite{BSqueeze} in the NLS setting. We plan to address this
topic in a forthcoming paper.

\subsection{Outline}

Sections 2 and 3 describe the multilinear correction technique
which generates modified energies. Section 4 establishes useful
pointwise upper bounds on certain multipliers arising in the multilinear
correction procedure. These upper bounds are combined with a quintilinear
estimate, in the $\R$ setting, to prove the bulk of \eqref{betteridea}
in Section 5. Section 6 contains the variant local well-posedness 
result and the proof of global well-posedness for \eqref{KdV} in 
$H^s ( \R ),~s > -\frac{3}{4}.$ We next consider the periodic initial
value problem \eqref{TKdV} with period $\lam$. Section 7 extends the
local well-posedness theory for \eqref{TKdV} to the $\lam$-periodic 
setting. Section 8 proves global well-posedness of \eqref{TKdV} in
$H^s ( \T ), ~s \geq - \half$. The last section exploits Miura's
transform to prove the corresponding global well-posedness results 
for the focussing and defocussing modified KdV equations.

\subsection{Notation}
We will use $c,C$ to denote various time independent constants, usually 
depending only upon $s$. In
case a constant depends upon other quantities, we will try to make that
explicit. We use $A \lesssim B$ to denote an estimate of the form $A \leq
C B$. Similarly, we will write $A \thicksim B$ to mean $A \lesssim B$
and $B \lesssim A$. To avoid an issue involving a logarithm, we depart
from standard practice and write 
$\langle k \rangle = 2 + |k|.$
The notation $a +$ denotes $a+ \epsilon $ for an arbitrarily small $\e$. 
Similarly, $a-$ denotes $a - \e$. 
We will make frequent use of the two-parameter spaces $X_{s,b} (\R \times \R)$
with norm 
\begin{equation}
  \label{Xsbdefined}
  \X s b u = {{\left( \int \int {{\langle k \rangle k}^{2s}} {{\langle \tau - \xi^3 \rangle}^{2b}}
{{| {\widehat{u}} ( \xi, \tau ) |}^2} d\xi d\tau \right)}^\half}.
\end{equation}
For any time interval $I$, we define the restricted spaces $X_{s,b} (R \times
I)$ by the norm
\begin{equation*}
  {{\| u \|}_{X_{s,b} ( \R \times I )}} = \inf \{ \X s b U : U|_{\R \times I}
= u \}.
\end{equation*}
We will systematically ignore constants involving $\pi$ in 
the Fourier transform, except in Section 7.
Other notation is introduced during the developments that follow.

\section{Multilinear forms}

In this section, we introduce notation for describing certain multilinear
operators, see for example \cite{CMWavelets}, \cite{CMOndelettes}. 
Bilinear versions of these operators will generate a sequence of
almost conserved quantities involving higher order multilinear corrections.

\begin{definition}
A {\it{k-multiplier}} is a function $m: \R^k \longmapsto \C$. A $k$-multiplier
is {\it{symmetric}} if $m( \xi) = m( \sigma( \xi ))$ for all $\sigma
\in S_n$, the group of all permutations on $n$ objects. The 
{\it{symmetrization}} of a $k$-multiplier $m$ is the multiplier
\begin{equation}
\label{symdefined}
[m]_{sym} ( \xi ) = \frac{1}{n!} \sum_{\sigma \in S_n } m ( \sigma ( \xi ) ).
\end{equation}
\end{definition}
The domain of $m$ is $\R^k$, however, we will only be interested in $m$
on the hyperplane $\xi_1 + \dots + \xi_k = 0$.

\begin{definition}
A $k$-multiplier generates a {\it{k-linear functional or k-form}} acting on $k$
functions $u_1, \dots, u_k$,
\begin{equation}
\label{Lambdakdefined}
\Lambda_k ( m; u_1, \dots, u_k) = \int\limits_{\xi_1 + \dots + \xi_k = 0}
m( \xi_1 , \dots , \xi_k ) {\widehat{u_1}} ( \xi_1 ) \dots
{\widehat{u_k}} ( \xi_k ) .
\end{equation}
We will often apply $\Lam_k $ to $k$ copies of the same function $u$
in which case the dependence upon $u$ may be suppressed in the notation:
$\Lam_k (m; u, \dots, u)$ may simply be written $\Lam_k (m)$.

If $m$ is symmetric then $\Lambda_k ( m )$ is a {\it{symmetric}}
$k$-linear functional.
\end{definition}

As an example, suppose that $u$ is a $\R$-valued function. We
calculate ${{\| u \|}_{L^2}^2} = \int \widehat{u} (\xi )
{\overline{\widehat{u}}} ( \xi ) d\xi = \int\limits_{\xi_1 + \xi_2 = 0}
\widehat{u} (\xi_1 ) \widehat{u} (\xi_2 ) = \Lam_2 (1).$

The time derivative of a symmetric $k$-linear functional can be calculated
explicitly if we assume that the function $u$ satisfies a particular PDE.
The following statement may be directly verified by using the KdV
equation.

\begin{proposition}
Suppose $u$ satisfies the KdV equation \eqref{KdV} and that $m$ is a
symmetric $k$-multiplier. Then
\begin{equation}
\label{timederivative}
\frac{d}{dt} \Lambda_k ( m) = \Lambda_k ( m \alpha_k ) - i \frac{k}{2}
\Lambda_{k+1} \left( m( \xi_1 , \dots, \xi_{k-1} , \xi_k + \xi_{k+1} \right) 
\{ \xi_k + \xi_{k+1} \}).
\end{equation}
where 
\begin{equation}
\label{alphak}
\alpha_k = i  ( \xi_1^3 + \dots + \xi_k^3 ).
\end{equation}
\end{proposition}
Note that the second term in \eqref{timederivative} may be symmetrized.

\section{Modified energies}

Let $m: \R \longmapsto \R$ be an arbitrary even $\R$-valued 1-multiplier
and define the associated operator by
\begin{equation}
{{\widehat{If}}} ( \xi ) = m( \xi) {\widehat{f}} ( \xi ).
\label{Ioperator}
\end{equation}
We define the {\it{modified energy}} $E_I^2 (t)$ by
\begin{equation*}
E_I^2 ( t) = {{\| Iu (t) \|}_{L^2}^2}.
\end{equation*}
The name ``modified energy'' is in part justified since 
in case $m=1, ~E^2_I (t) = {{\| u(t) \|}_{L^2}^2}.$ We will show
later that for $m$ of a particular form, certain modified
energies enjoy an almost conservation property.
By Plancherel and the fact that $m$ and $u$ are $\R$-valued,
\begin{equation*}
E_I^2 ( t) = \Lambda_2 ( m( \xi_1 ) m( \xi_2 ) ).
\end{equation*}
Using \eqref{timederivative}, we have
\begin{equation}
\label{firsttimed}
\frac{d}{dt} E^2_I ( t) = \Lambda_2 ( m( \xi_1 ) m( \xi_2 ) \alpha_2 )
- i \Lambda_3 ( m( \xi_1 ) m( \xi_2 + \xi_3 ) \{ \xi_2 + \xi_3 \} ).
\end{equation}
The first term vanishes. We symmetrize the remaining
term to get
\begin{equation*}
\frac{d}{dt} E^2_I ( t) = \Lambda_3 ( -i [ m(\xi_1 ) m( \xi_2 + \xi_3 )
(\xi_2 + \xi_3 ) ]_{sym} ).
\end{equation*}
Note that the time derivative of $E^2_I (t) $ is a 3-linear expression.
Let us denote
\begin{equation}
\label{Mthreedefined}
M_3 ( \xi_1 , \xi_2 , \xi_3 ) = -i [ m(\xi_1 ) m( \xi_2 + \xi_3 )
\{ \xi_2 + \xi_3 \} ]_{sym} .
\end{equation}
Observe that if $m =1 $, the symmetrization results in $M_3 = c (\xi_1 + \xi_2
+ \xi_3 )$. This reproduces the Fourier proof of $L^2$ mass conservation
from the introduction.

Form the new modified energy
\begin{equation*}
E^3_I (t) = E^2_I ( t) + \Lambda_3 ( \sigma_3 )
\end{equation*}
where the symmetric 3-multiplier $\sigma_3$ will be chosen momentarily
to achieve a cancellation. Applying \eqref{timederivative} gives
\begin{equation}
\label{secondtimed}
\frac{d}{dt} E_I^3 (t) = \Lambda_3 ( M_3 ) + \Lambda_3 ( \sigma_3 \alpha_3 )
+ \Lambda_4 \left( - i \frac{3}{2} \sigma_3 ( \xi_1 , \xi_2,
\xi_3 + \xi_4 ) \{ \xi_3 + \xi_4 \}  \right) .
\end{equation}
We choose
\begin{equation}
\label{sigmathreedefined}
\sigma_3 = - \frac{M_3}{\alpha_3}
\end{equation}
to force the two $\Lambda_3 $ terms in \eqref{secondtimed} to cancel. With
this choice, the time derivative of $E^3_I (t) $ is a 4-linear expression
$\Lam_4 ( M_4 )$ where
\begin{equation} 
\label{Mfourdefined}
M_4 ( \xi_1 , \xi_2 , \xi_3 , \xi_4 ) = - i \frac{3}{2}
[ \sigma_3 ( \xi_1 , \xi_2 , \xi_3 + \xi_4 ) \{ \xi_3 + \xi_4 \} ]_{sym} .
\end{equation}
Upon defining
\begin{equation*}
E_I^4 (t) = E^3_I (t) + \Lambda_4 ( \sigma_4)
\end{equation*}
with 
\begin{equation}
\label{sigmafourdefined}
\sigma_4 = - \frac{M_4}{\alpha_4}
\end{equation}
we obtain
\begin{equation}
\label{Efourinc}
\frac{d}{dt} E^4_I (t) = \Lambda_5 ( M_5 )
\end{equation}
where
\begin{equation}
\label{Mfivedefined}
M_5 ( \xi_1 , \dots, \xi_5 ) = -2i [
\sigma_4 ( \xi_1, \xi_2 , \xi_3 , \xi_4 + \xi_5 ) \{ \xi_4 + \xi_5 \} ]_{sym}.
\end{equation}
This process can clearly be iterated to generate $E^n_I$ satisfying
$\frac{d}{dt} E^n_I (t) = \Lambda_{n+1} ( M_{n+1} )$, $n = 2, 3 , \dots.$
These higher degree corrections to the modified energy $E_I^2$ may be of
relevance in studying various qualitative aspects of the KdV evolution.
However, for the purpose of showing GWP in $H^s ( \R )$ down to 
$s > - \frac{3}{4}$ and in $H^s (\T )$ down to $s \geq - \half$, we 
will see that almost conservation of $E^4_I (t)$
suffices.

The modified energy construction process is illustrated in the case of the Diriclet energy
\begin{equation*}
  E^2_D (t) = {{\| \partial_x u \|}_{L^2_x}^2} = \Lam_2 ( (i\xi_1 ) (i \xi_2)).
\end{equation*}
Define $E^3_D (t) = E^2_D (t) + \Lam_3 ( \sigma_3 )$, 
and use \eqref{timederivative} to see
\begin{equation*}
  \partial_t E^3_D (t) = \Lam_3 ( [ i(\xi_1 + \xi_2 ) i \xi_3 \{ \xi_1 + \xi_2 \}
]_{sym} ) + \Lam_3 ( \sigma_3 \alpha_3 ) + \Lambda_4 (M_4),
\end{equation*}
where $M_4$ is explicitly obtained from $\sigma_3$. Noting that
$i (\xi_1 + \xi_2 ) i \xi_3 \{ \xi_1 + \xi_2 \} = - \xi_3^3 $ on the set 
$\xi_1 + \xi_2 + \xi_3 =0$, we know that
\begin{equation*}
  \partial_t E^3_D 
(t) = \Lam_3 (-\frac{1}{3} \alpha_3 ) + \Lam_3 (\sigma_3 \alpha_3 ) + \Lam_4 ( M_4 ).
\end{equation*}
The choice of $\sigma_3 = \frac{1}{3}$ results in a cancellation of the
$\Lam_3$ terms and $M_4 = [\{ \xi_1 + \xi_2 \}]_{sym} = \xi_1 + \xi_2 + \xi_3 
+ \xi_4$ so $M_4 = 0$.

Therefore, $E^3_D 
(t) = \Lam_2 ( (i \xi_1 ) (i \xi_2 )) + \Lam_3 ( \frac{1}{3} )$
is an exactly conserved quantity. The modified energy construction applied
to the Dirichlet energy led us to the Hamiltonian for KdV. Applying the
construction to higher order derivatives in $L^2$ will similarly 
lead to the higher
conservation laws of KdV.

\section{Pointwise multiplier bounds}

This section presents a detailed analysis of the multipliers $M_3, ~ M_4, ~
M_5 $ which were introduced in the iteration process of the previous section.
The analysis identifies cancellations resulting in pointwise upper bounds on 
these multipliers depending upon the relative sizes of the multiplier's
arguments. These bounds are applied to prove an almost conservation
property in the next section. 
We begin by recording some arithmetic and calculus facts. 

\subsection{Arithmetic and calculus facts}

The following arithmetic facts may be easily verified:
\begin{equation}
\xi_1 + \xi_2 + \xi_3 = 0 \implies \alpha_3 = \xi_1^3 + \xi_2^3
+ \xi_3^3 = 3 \xi_1 \xi_2 \xi_3.
\label{factthree}
\end{equation}
\begin{equation}
\xi_1 + \xi_2 + \xi_3 + \xi_4 = 0 \implies \alpha_4 = \xi_1^3 + \xi_2^3
+ \xi_3^3 + \xi_4^3 = 3 (\xi_1 + \xi_2) ( \xi_1 + \xi_3 ) (\xi_1 + \xi_4).
\label{factfour}
\end{equation}
A related observation for the circle was exploited by C. Fefferman 
\cite{FeffSphericalSummation} and by Carleson and
Sjolin \cite{CarSj72} for curves with nonzero curvature. 
These properties were also observed by
Rosales \cite{RosalesThesis} and \eqref{factthree} was used by Bourgain in
\cite{B1}.

\begin{definition} Let $a$ and $b$ be smooth functions of the real
variable $\xi$. We say that {\it{ $a$ is controlled by $b$}} if $b$ is 
non-negative and satisfies $b(\xi) \thicksim b(\xi' )$ 
for $|\xi | \thicksim |\xi' |$
and
\begin{eqnarray*}
a(\xi) &=& O ( b(\xi ) ), \\
a'(\xi ) &=& O \left( \frac{b(\xi)}{|\xi |} \right), \\
a''(\xi ) &=& O \left( \frac{b(\xi )}{|\xi|^2} \right), \\
\end{eqnarray*}
for all non-zero $\xi$.
\end{definition}
 
With this notion, we can state the following forms of the mean value theorem.

\begin{lemma} 
\label{onemvt}
If $a$ is controlled by $b$ and $|\eta | \ll |\xi |$, then
\begin{equation}
\label{mvt}
a( \xi + \eta ) - a( \xi ) = O \left( |\eta| \frac{b(\xi )}{|\xi |} \right).
\end{equation}
\end{lemma}

\begin{lemma}
\label{twomvt}
If $a$ is controlled by $b$ and $|\eta |,~ |\lam | \ll |\xi |$ then
\begin{equation}
\label{doublemvt}
a(\xi + \eta + \lam ) - a( \xi + \eta ) - a( \xi + \lam ) - a ( \xi )
= O \left( |\eta | |\lam | \frac{b( \xi )}{|\xi |^2 } \right).
\end{equation}
\end{lemma}
We will sometimes refer to our use of \eqref{doublemvt} as applying
the {\it{double mean value theorem}}. 
\subsection{$M_3$ bound}

The multiplier $M_3$ was defined in \eqref{Mthreedefined}. In this
section, we will generally be considering an arbitrary even $\R$-valued
1-multiplier $m$. We will specialize to the situation when $m$
is of the form \eqref{particularm} below. Recalling 
that $\xi_1 + \xi_2 + \xi_3 =0$ and that 
$m$ is even allows us to reexpress \eqref{Mthreedefined} as
\begin{equation}
\label{Mthreeviam}
M_3 ( \xi_1 , \xi_2 , \xi_3 ) = -i [ m^2 (\xi_1  ) \xi_1 ]_{sym}
= -\frac{i}{3} [ m^2 ( \xi_1 ) \xi_1 + m^2 ( \xi_2 ) \xi_2 + 
 m^2 ( \xi_3 ) \xi_3 ].
\end{equation}

\begin{lemma}
\label{Mthreelemma}
If $m$ is even, $\R$-valued and 
$m^2$ is controlled by itself then, on the set $\xi_1 + \xi_2 + \xi_3 =0,~
|\xi_i | \thicksim N_i$ (dyadic),
\begin{equation}
\label{Mthreebound}
|M_3 ( \xi_1 , \xi_2 , \xi_3 ) | \lesssim \max ( m^2 ( \xi_1 ), 
m^2 ( \xi_2 ) , m^2 ( \xi_3 )) \min ( N_1 , N_2 , N_3 ).
\end{equation}
\end{lemma}

\begin{proof} {%of Mthreelemma
Symmetry allows us to assume $N_1 = N_2 \geq N_3$. In case $N_3 \ll N_1$, 
the claimed estimate
is equivalent to showing
\begin{equation*}
m^2 ( \xi_1 ) \xi_1 - m^2 ( \xi_1 + \xi_3 ) ( \xi_1 + \xi_3 )
+ m^2 (\xi_3 ) \xi_3 \leq \max ( m^2 (\xi_1 ) , m^2 ( \xi_3 )) N_3.
\end{equation*}
But this easily follows when we rewrite the left-side as $(m^2( \xi_1 ) -
m^2 ( \xi_1 + \xi_3 ) ) \xi_1 - m^2 ( \xi_1 + \xi_3 ) \xi_3 + m^2 (\xi_3 ) 
\xi_3  $ and use \eqref{mvt}.} In case $N_3 \thicksim N_2$, \eqref{Mthreebound}
may be directly verified.
\end{proof}
In the particular case when the multiplier $m( \xi )$ is smooth, monotone,
and of the form
\begin{equation}
\label{particularm}
m(\xi) = \left\{ \begin{matrix} 1, & |\xi | < N, \\ 
                              N^{-s} |\xi|^s, & |\xi | > 2N,
\end{matrix}
\right.  
\end{equation}
we have 
\begin{equation}
\label{particularthree}
|M_3 ( \xi_1 , \xi_2 , \xi_3 ) |\leq \min( N_1, N_2 , N_3) .
\end{equation}

\subsection{$M_4$ bound.}

This subsection establishes the following pointwise upper bound on the
multiplier $M_4$.

\begin{lemma}
  \label{MFourlemma}
Assume $m$ is of the form \eqref{particularm}.
In the region where $|\xi_ i | \thicksim N_i , ~|\xi_j + \xi_k | \thicksim 
N_{jk}$ for $N_i , N_{jk} $ dyadic,
\begin{equation}
\label{MFourEst}
|M_4 ( \xi_1 , \xi_2 , \xi_3 , \xi_4 )| \lesssim \frac{ |\alpha_4 |  
~ m^2 ( \min
( N_i , N_{jk} ) ) }{(N+N_1) (N+N_2) (N+N_3) (N+N_4)} .
\end{equation}
\end{lemma}

We begin by deriving two explicit representations of $M_4$ in terms of $m$.
These identities are then analyzed in cases to prove \eqref{MFourEst}.

Recall that, 
\begin{equation}
  \label{Mfour}
  M_4 ( \xi_1 , \xi_2 , \xi_3 , \xi_4 ) = c [ \sigma_3 ( \xi_1 , \xi_2, \xi_3
+ \xi_4 ) (\xi_3 + \xi_4 ) ]_{sym},
\end{equation}
where $\sigma_3 = - \frac{M_3}{\alpha_3}$ and
\begin{equation}
  \label{Mthree}
  M_3 ( x_1 , x_2 , x_3 ) = -i [ m(x_1 ) m( x_2 + x_3 )(x_2 + x_3 ) ]_{sym}
= - \frac{i}{3}[m^2 (x_1 ) x_1 + m^2 (x_2) x_2 + m^2 (x_3 ) x_3] , 
\end{equation}
and $\alpha_3 (x_1 , x_2 , x_3 ) = x_1^3 + x_2^3 + x_3^3 = 3 x_1 x_2 x_3 $.
We shall ignore the irrelevant constant in \eqref{Mfour}.
Therefore,
\begin{equation}
\label{Mfromm}
   M_4 ( \xi_1 , \xi_2 , \xi_3 , \xi_4 ) = -\half \left[ \frac{ 
m^2 (\xi_1 ) \xi_1  + 
m^2 (\xi_2) \xi_2 + m^2 (\xi_3 + \xi_4 )(\xi_3 + \xi_4 ) }{
3 \xi_1 \xi_2 } \right]_{sym}  
\end{equation}
\begin{equation*}
= -\half \left[ \frac{ 
2 m^2 (\xi_1 ) \xi_1  + m^2 (\xi_3 + \xi_4 )(\xi_3 + \xi_4 ) }{
3 \xi_1 \xi_2 } \right]_{sym} .
\end{equation*}
Recall also from \eqref{factfour} that
\begin{equation}
\label{alphafour}
\alpha_4 ( \xi_1 , \xi_2 , \xi_3 , \xi_4 ) =
\xi_1^3 + \xi_2^3 + \xi_3^3 + \xi_4^3 = 3 (\xi_1 \xi_2 \xi_3
+ \xi_1 \xi_2 \xi_4 + \xi_1 \xi_3 \xi_4 + \xi_2 \xi_3 \xi_4 )
\end{equation}
\begin{equation*}
    = 3 ( \xi_1 + \xi_2 ) (\xi_1 + \xi_3 ) (\xi_1 + \xi_4 ).
\end{equation*}

We can now rewrite the first term in \eqref{Mfromm}
\begin{equation}
   \left[ \frac{ 2 m^2 ( \xi_1 ) \xi_1 \xi_3 \xi_4 }{3 
\xi_1 \xi_2 \xi_3 \xi_4 }
\right]_{sym}
= \frac{2}{9}
\left[ \frac{m^2(\xi_1 ) (\xi_1 \xi_2 \xi_3
+ \xi_1 \xi_2 \xi_4 + \xi_1 \xi_3 \xi_4 + \xi_2 \xi_3 \xi_4  - 
\xi_2 \xi_3 \xi_4) }{\xionetofour} \right]_{sym} 
 \label{Aterm}
\end{equation}
\begin{equation*}
= \frac{1}{54} \left[ m^2 ( \xi_1 ) + m^2 ( \xi_2 ) + 
m^2 ( \xi_3 ) + m^2 ( \xi_4 ) \right] \frac{\alpha_4 }{\xionetofour}
- \frac{1}{18} \left[ \frac{m^2(\xi_1 )}{\xi_1 } + \frac{m^2(\xi_2 )}{\xi_2 }
+ \frac{m^2(\xi_3 )}{\xi_3 }+  \frac{m^2(\xi_4 )}{\xi_4 } \right].
\end{equation*} 
The second term in \eqref{Mfromm} is rewritten, using $\xi_1 + \xi_2 + 
\xi_3 + \xi_4 =0$, and the fact the $m$ is even,
\begin{equation*}
\left[ \frac{ - m^2 (\xi_1 + \xi_2 ) (\xi_1+ \xi_2 ) \xi_3 
\xi_4 }{3 \xionetofour}  \right]_{sym} = 
- \frac{1}{18} \left\{ \frac{ m^2( \xi_1 + \xi_2 ) (\xi_1 \xi_3 \xi_4 + \xi_2
\xi_3 \xi_4) + m^2 ( \xi_3 + \xi_4 )  (\xi_1 \xi_2 \xi_3 + \xi_1
\xi_3 \xi_4)}{\xionetofour}  \right.
\end{equation*}
\begin{equation*}
+ \frac{ m^2( \xi_1 + \xi_3 ) (\xi_1 \xi_2 \xi_4 + \xi_2
\xi_3 \xi_4) + m^2 ( \xi_2 + \xi_4 )  (\xi_1 \xi_2 \xi_3 + \xi_1
\xi_3 \xi_4)}{\xionetofour}
\end{equation*}
\begin{equation*} 
\left.
+ \frac{ m^2( \xi_1 + \xi_4 ) (\xi_1 \xi_2 \xi_3 + \xi_2
\xi_3 \xi_4) + m^2 ( \xi_2 + \xi_3 )  (\xi_1 \xi_2 \xi_4 + \xi_1
\xi_3 \xi_4)}{\xionetofour}  \right\}
\label{Bterm}
\end{equation*}
\begin{equation}
\label{Brewr}
= - \frac{1}{54} \frac{\alpha_4 }{\xionetofour} \left[
m^2 ( \xi_1 + \xi_2 ) + m^2 ( \xi_1 + \xi_3 ) + m^2 ( \xi_1 + \xi_4 ) \right].
\end{equation}
 
We record two identities for $M_4$.

\begin{lemma} If $m$ is even and $\R$-valued, the following two identities
for $M_4$ are valid:
\label{MFourid}
\begin{equation}
  \label{oneseven}
  M_4 (\xi_1 , \xi_2 ,\xi_3 , \xi_4 )
= - \frac{1}{108}  \frac{\alpha_4 }{\xionetofour} \left[
m^2 (\xi_1 ) + m^2 (\xi_2 )+ m^2 (\xi_3 ) + m^2 (\xi_4 ) \right.
\end{equation}
\begin{equation*}
\left.
- m^2 (\xi_1 + \xi_2 ) -  m^2 (\xi_1 + \xi_3 )  - m^2 (\xi_1 + \xi_4 ) \right]
\end{equation*}
\begin{equation*}
  + \frac{1}{36} \left\{ \frac{m^2 (\xi_1 ) }{\xi_1 } 
+ \frac{m^2 (\xi_2 ) }{\xi_2 } +
\frac{m^2 (\xi_3 ) }{\xi_3 } + \frac{m^2 (\xi_4 ) }{\xi_4 } \right\} :=I + II.
\end{equation*}
\begin{equation}
  \label{oneten}
  M_4 ( \xi_1 , \xi_2 , \xi_3 , \xi_4 )
= -\frac{1}{36} \frac{1}{\xi_1 \xi_2 \xi_3 \xi_4 } \times
\end{equation}
\begin{eqnarray*}
& \{ & \xi_1 \xi_2 \xi_3 [ m^2 ( \xi_1 ) + m^2 ( \xi_2 ) + m^2 ( \xi_3 ) 
- m^2 ( \xi_1 + \xi_2 ) - m^2 ( \xi_1 + \xi_3 ) - m^2 ( \xi_1 + \xi_4 ) ]  \\
&+& \xi_1 \xi_2 \xi_4 [ m^2 ( \xi_1 ) + m^2 ( \xi_2 ) + m^2 ( \xi_4 ) 
- m^2 ( \xi_1 + \xi_2 ) - m^2 ( \xi_1 + \xi_3 ) - m^2 ( \xi_1 + \xi_4 ) ] \\
&+&
 \xi_1 \xi_3 \xi_4 [ m^2 ( \xi_1 ) + m^2 ( \xi_3 ) + m^2 ( \xi_4 ) 
- m^2 ( \xi_1 + \xi_2 ) - m^2 ( \xi_1 + \xi_3 ) - m^2 ( \xi_1 + \xi_4 ) ] \\
&+& \xi_2 \xi_3 \xi_4 [ m^2 ( \xi_2 ) + m^2 ( \xi_3 ) + m^2 ( \xi_4 ) 
- m^2 ( \xi_1 + \xi_2 ) - m^2 ( \xi_1 + \xi_3 ) - m^2 ( \xi_1 + \xi_4 ) ] \}.\\
\end{eqnarray*}
\end{lemma}
\begin{proof}
The identity \eqref{oneseven} was established above. The identity 
\eqref{oneten} follows from \eqref{oneseven} upon expanding $\alpha_4$
and writing the second term in \eqref{oneseven} on a common denominator. 
\end{proof}

\begin{proof} [Proof of Lemma \ref{MFourlemma}]
The proof consists of a case-by-case analysis pivoting on the relative sizes
of $N_i, N_{jk}$. Symmetry properties of $M_4$ permit us to assume that
$|\xi_1 | \geq |\xi_2 | \geq |\xi_3 | \geq |\xi_4 |$. Consequently, we
assume $N_1 \geq N_2 \geq N_3 \geq N_4.$ Since $m^2 ( \xi) =1$
for $|\xi| < \frac{N}{2}$, a glance at \eqref{Mfromm} shows that $M_4$
vanishes when $|\xi_1 | < \frac{N}{4}.$ We may therefore assume that $|\xi_1 |
\gtrsim N$. Since $\xi_1 + \xi_2 + \xi_3 + \xi_4 = 0$, we must also
have $|\xi_2 | \gtrsim N$.

From \eqref{alphafour}, we know that we can replace $\alpha_4$ on the
right side of \eqref{MFourEst} by $N_{12} N_{13} N_{14}$. Suppose $N_{12}
< \frac{N_1}{2}, ~N_{13}< \frac{N_1}{2}, ~N_{14}< \frac{N_1}{2}.$ Then,
$\xi_1 \thicksim - \xi_2, ~\xi_1 \thicksim - \xi_3$ and $\xi_1 \thicksim
- \xi_4 $ so $\xi_1 + \xi_2 + \xi_3 + \xi_4 \thicksim -2\xi_1 \neq 0$.
Thus, at least one of $N_{12} , ~N_{13} , ~N_{14}$ must be at least of size 
comparable to $N_1$. The right side of
\eqref{MFourEst} may be reexpressed
\begin{equation}
\label{reexprhs}
\frac{N_{12} N_{13} {N_{14}} 
m^2 ( \min ( N_i , N_{jk} )) }{ N_1^2 (N+ N_3 ) (N+ N_4 ) }.
\end{equation}

{\bf{Case 1.}} $| N_4 | \gtrsim \frac{N}{2}$. \newline
Term $I$ in \eqref{oneseven} is bounded by 
$\frac{N_{12} N_{13} N_{14} }{ N_1^2 
N_3 N_4 } m^2 ( \min (N_i , N_{jk} )) $, 
and therefore, after cancelling $\max( N_{12}, N_{13} , N_{14})$
with one of the $N_1$, satisfies \eqref{MFourEst}. Term $II$ is treated next.
In case $N_{12}, N_{13} , N_{14} \gtrsim N_1$, \eqref{reexprhs} is an upper
bound of $\frac{N_1}{N_3 N_4} m^2 (N_4 ) \geq \frac{m^2 ( N_4 ) }{N_4 }$
and the triangle inequality gives $|II| \lesssim \frac{m^2 ( N_4 )}{N_4 }$
since $\frac{m^2 ( \cdot )}{( \cdot )}$ is a decreasing function. If $N_{12}
\gtrsim N_1, ~ N_{13} \ll N_1 $ and $N_{14} \gtrsim N_1 $, we rewrite
\begin{equation*}
|II| \thicksim \left\{ \frac{m^2 (\xi_1 ) }{\xi_1 } + \frac{
m^2 ( - \xi_1 + ( \xi_1 + \xi_3 ))}{(-\xi_1 + (\xi_1 + \xi_3 ))} +
\frac{m^2 ( \xi_2 ) }{\xi_2 } +  \frac{
m^2 ( - \xi_2 + ( \xi_2 + \xi_4 ))}{(-\xi_2 + (\xi_2 + \xi_4 ))} \right\}.
\end{equation*}

Applying the mean value theorem and using $\xi_1 + \xi_2 + \xi_3 + \xi_4 = 0$ 
gives $|II| \lesssim {{\left( \frac{m^2 
( \wt{\xi_1 } )}{ \wt{\xi_1 } } \right)}'} ( \xi_1 + \xi_3 ) \lesssim
\frac{N_{13}}{N_1^2 } m^2 (N_1) $ since $\wt{\xi_1 } = \xi_1 + O(N_{13} )$ 
and $N_{13} \ll N_1$, so this subcase is fine. If $N_{12} \ll N_1 , 
N_{13} \ll N_1 $ and $N_{14} \gtrsim N_1 $, the double mean value 
theorem \eqref{doublemvt} applied to term $II$ gives
the bound
\begin{equation*}
|II| \thicksim {{\left( 
\frac{m^2 ({{ \xi_1 } })}{ {{\xi_1 }^3} } \right)}''} (\xi_1 + \xi_2 )
(\xi_1 + \xi_3 )
\end{equation*}
Our assumptions
on $N_{12}, ~ N_{13}$ give
the bound $|II| \lesssim \frac{N_{12} N_{13} }{N_1^3} m^2 ( N_1 )$ which is 
smaller than \eqref{reexprhs}.

The remaining subcases have either precisely one element of the set $\{ N_{12},
N_{13}, N_{14} \}$ much smaller than $N_1$ or precisely two elements much
smaller than $N_1$. In the case of just one small $N_{1j}$, we apply the
mean value theorem as above. When there are two small $N_{1j}$, we apply
the double mean value theorem as above.

{\bf{Case 2.}} $| N_4 | \ll \frac{N}{2}. $ \newline
Certainly, $m^2 ( \min ( N_i , N_{jk} )) = 1$ in this region. It is
not possible for both $N_{12} < \frac{N_1}{4}$ and $N_{13} < \frac{N_1}{4}$
in this region. Indeed, we find then that $\xi_1 \thicksim - \xi_2$ and
$\xi_1 \thicksim - \xi_3 $ which with $\xi_1 + \xi_2 + \xi_3 +\xi_4 = 0$
implies $\xi_4 \thicksim \xi_1 $ but $|\xi_4 | \ll \frac{N}{2}$ while
$|\xi_1 | \thicksim N_1 \gtrsim N$. 
We need to show $M_4 \leq \frac{N_{12} N_{13} }{N_1 (N+N_3 ) N }$.
\newline
{\bf{Case 2.A.}} $\frac{N_1}{4} > 
N_{12} \gtrsim \frac{N}{2},~ N_{13} \thicksim N_1 $. \newline
Since $N_4 \ll \frac{N}{2}$ and $\xi_1 + \xi_2 + \xi_3 + \xi_4 = 0$, we
must have $N_{12} \thicksim N_3 $. So $N+ N_3 \thicksim N_3 $ and our
goal is to show $M_4 \lesssim \frac{N_{12}}{N_3 N } \thicksim \frac{1}{N}.$
The last three terms in \eqref{oneten} are all $O(\frac{1}{N })$ which
is fine. The first term in \eqref{oneten} is
\begin{equation*}
\frac{1}{18 \xi_4 } ( m^2 ( \xi_1 ) + m^2 ( \xi_2 ) + m^2 ( \xi_3 )
- m^2 ( \xi_1 + \xi_2 ) - m^2 ( \xi_1 + \xi_3 ) - m^2 ( \xi_1 + \xi_4 ) ) .
\end{equation*} 
Replacing $\xi_1 + \xi_2$ by  $-(\xi_3 + \xi_4)$ and $\xi_1 + \xi_3 $
by $-(\xi_2 + \xi_4 )$, we identify three differences poised for the mean
value theorem. We find this term equals
\begin{equation*}
\frac{1}{18 \xi_4 } [ (m^2 ( \wt{\xi_1 } ) )' + (m^2 ( \wt{\xi_2 } ) )' +  
(m^2 ( \wt{\xi_3 } ) )' ] \xi_4
\end{equation*}
with $\wt{\xi_i} = \xi_i + O(N_4 ) $ for $i = 1,2,3$ 
so $|\wt{\xi_i }| \thicksim N_i $. This
expression is also $O(\frac{1}{N } )$. 
\newline
{\bf{Case 2.B.}} $ N_{12} \ll \frac{N}{2},~ N_{13} \thicksim N_1$. \newline
Since $N_{12} = N_{34}$ and $N_4 \ll \frac{N}{2}$, we must have $N_3 \ll
\frac{N}{2}$. We have $N_{13} \thicksim N_1$ and $N_{14} \thicksim N_1$
here so our desired upper bound is $\frac{N_{12}}{N^2}$. We recall
\eqref{oneseven} and evaluate $m^2 $ when we can to find
\begin{equation}
\label{twob}
M_4 ( \xi_1 , \xi_2 , \xi_3 , \xi_4 ) = \frac{\alpha_4 }{54 \xi_1
\xi_2 \xi_3 \xi_4 } (m^2 ( \xi_1 ) + m^2 ( \xi_2 ) +2 -1 - m^2 (\xi_1
+ \xi_3 ) - m^2 ( \xi_1 + \xi_4 ) )
\end{equation}
\begin{equation*}
- \frac{1}{18} \left( \frac{m^2 ( \xi_1 )}{\xi_1 } +  
\frac{m^2 ( \xi_2 )}{\xi_2 } + \frac{\xi_3 + \xi_4 }{\xi_3 \xi_4 } \right).
\end{equation*}
The last term is dangerous so we isolate a piece of the first term to
cancel it away. Expanding $\alpha_4 = 3 (\xi_1 + \xi_2 ) (\xi_1 + \xi_3 )
(\xi_1 + \xi_4 ) $, we see that
\begin{eqnarray*}
\frac{\alpha_4}{54 \xi_1 \xi_2 \xi_3 \xi_4 } &=& \frac{1}{18}
\frac{ (\xi_3 + \xi_4 ) }{\xi_3 \xi_4 } 
\frac{(\xi_2 + \xi_4 )(\xi_1 + \xi_4 )}{\xi_1 \xi_2 } \\
&=& \frac{1}{18} \frac{ (\xi_3 + \xi_4 ) }{\xi_3 \xi_4 } 
\left( 1 + \frac{\xi_4 ( \xi_1 + \xi_2 + \xi_4 )}{\xi_1 \xi_2 } \right) \\
&=& \frac{1}{18}  \frac{ (\xi_3 + \xi_4 ) }{\xi_3 \xi_4 }
\left( 1 - \frac{\xi_4 \xi_3 }{\xi_1 \xi_2 } \right). 
\end{eqnarray*}
The first piece cancels with $- \frac{1}{18} \frac{ \xi_3 + \xi_4 }{\xi_3 
\xi_4 }$ in \eqref{twob} and the second piece is of size 
$\frac{N_{12}}{N_1^2}$ which is fine. It remains to control
\begin{equation}
\label{herestep}
\frac{\alpha_4}{54 \xi_1 \xi_2 \xi_3 \xi_4 }
\left( m^2 ( \xi_1 ) + m^2 (\xi_2 ) - m^2 ( \xi_1 + \xi_3 ) - m^2 ( \xi_1
+ \xi_4 ) \right) - \frac{1}{18} \left( \frac{m^2 ( \xi_1 )}{\xi_1 } +
\frac{ m^2 ( \xi_2 ) }{\xi_2 } \right),
\end{equation}
by $\frac{N_{12}}{N^2} .$ Expand $\alpha_4 $ using \eqref{alphafour}
to rewrite this expression as
\begin{equation}
\label{almostdone}
\frac{3 (\xi_1 \xi_2 \xi_3 + \xi_1 \xi_2 \xi_4 ) }{54 \xi_1 \xi_2 \xi_3 
\xi_4 } \left( m^2 ( \xi_1 ) + m^2 ( \xi_2 ) - m^2 ( \xi_1 + \xi_3 )
- m^2 ( \xi_1 + \xi_4 ) \right) 
\end{equation}
\begin{equation*}
- \frac{3 ( \xi_1 \xi_3 \xi_4 + \xi_2 \xi_3 \xi_4 ) }{54 \xi_1 \xi_2 
\xi_3 \xi_4 } ( m^2 ( \xi_1 + \xi_3 ) + m^2 ( \xi_1 + \xi_4 ) )
+ \frac{1}{18 \xi_1 \xi_2 } [ \xi_1 m^2 ( \xi_1 )+ \xi_2 m^2 ( \xi_2 )].
\end{equation*}
(The second term in \eqref{herestep} cancelled with part of the first.)
The second and third terms 
in \eqref{almostdone} are $O(\frac{N_{12} }{N^2})$ and may therefore be
ignored.
We rewrite the first term in \eqref{almostdone} using the fact that
$m^2$ is even as
\begin{equation*}
\frac{3 (\xi_1 \xi_2 \xi_3 + \xi_1 \xi_2 \xi_4 ) }{54 \xi_1 \xi_2 \xi_3
\xi_4 } \left( m^2 ( - \xi_1 ) + m^2 ( \xi_2 ) - m^2 ( - (\xi_1 + \xi_3 ))
- m^2 ( - (\xi_1 + \xi_4 )) \right).
\end{equation*}
Since $- \xi_1 + \xi_2 + (\xi_1 + \xi_3 ) + (\xi_1 + \xi_4 ) = 0$,
we can apply the double mean value theorem to obtain
\begin{equation*}
= \frac{3 ( \xi_1 \xi_2 \xi_3 + \xi_1 \xi_2 \xi_4 ) }{54 \xi_1 \xi_2 \xi_3 
\xi_4 } {{( m^2 ( - {\wt{\xi_1 }} ) )}''} \xi_3 \xi_4 
\end{equation*}
with $- \wt{\xi_1 } = - \xi_1 + O ( N_3 ) + O (N_4 ) \implies |- \wt{\xi_1 } |
\thicksim N_1$. Therefore, this term is bounded by
\begin{equation*}
\frac{ \xi_1 \xi_2 \xi_3 + \xi_1 \xi_2 \xi_4 }{\xi_1 \xi_2 \xi_3 \xi_4 }
\frac{m^2 ( - \wt{\xi_1 } )} { {{(- \wt{\xi_1 } )}^2} } \xi_3 \xi_4
= O \left( \frac{N_{12}}{N_1^2 } m^2 ( \wt{\xi_1 } ) \right),
\end{equation*}
which is smaller than $\frac{N_{12}}{N_1^2}$ as claimed.  
\newline
{\bf{Case 2.C.}} $\frac{N_1}{4} > N_{13} \gtrsim 
\frac{N}{2}, ~ N_{12} \thicksim
N_1$. \newline
This case follows from a modification of Case 2A.
\newline
{\bf{Case 2.D.}} $N_{13} \ll \frac{N}{2}, ~N_{12} \thicksim N_1 .$ \newline
This case does not occur because $N_{13} \thicksim N_{24}$ but $N_4$
is very small which forces $N_2$ to also be small which is a contradiction.
\end{proof} % of M_4 bound lemma
%%%%%%%%%%%%%%%%%%%%%%%%%%%%%%%%%%%%%%%%%%%%%%%%%%%%%%%%%%%%%%%%%%%%%%%%

\subsection{$M_5$ bound}

The multiplier $M_5$ was defined in \eqref{Mfivedefined}, with $\sigma_4
= - \frac{M_4}{\alpha_4}.$ Our work on $M_4$ above showed that $M_4$
vanishes whenever $\alpha_4$ vanishes so there is no denominator singularity
in $M_5$. Moreover, we have the following upper bound on $M_5$ in the
particular case when $m$ is of the form \eqref{particularm}.

\begin{lemma} If $m$ is of the form \eqref{particularm} then
\label{MFivelemma}
\begin{equation}
\label{MFivebound}
|M_5 ( \xi_1 , \dots, \xi_5) | \lesssim
\left[ \frac{ m^2 (N_{*45} ) ~N_{45} }{ (N+N_1 )
(N+N_2 )(N+N_3 ) (N+N_{45} ) }
\right]_{sym},
\end{equation}
where
\begin{equation*}
  N_{*45} = \min (N_1 , N_2 , N_3 , N_{45}, N_{12}, N_{13}, N_{23} ).
\end{equation*}
\end{lemma}
\begin{proof}
This follows directly from Lemma \ref{MFourlemma}. Note that 
$\xi_1 + \xi_4 + \xi_5 = - (\xi_2 + \xi_3 )$ allows
for the simplification in defining $N_{*45}$.
\end{proof}

\section{Quintilinear estimate on $\R$}

The $M_5$ upper bound contained in Lemma \ref{MFivelemma} and the local
wellposedness machinery \cite{KPVCPAM}, \cite{B1}, \cite{KPVBilin} are
applied to prove an almost conservation property of the modified energy
$E^4_I$. The almost conservation of $E^4_I$ is the key ingredient in
our proof of global wellposedness of the initial value problem for KdV
with rough initial data.

Let $X_{s,b}^\delta$ denote the Bourgain space \cite{B1} associated
to the cubic $\{ \tau = \xi^3 \}$ on the time interval $[0,\delta]$.
We begin with a quintilinear estimate.

\begin{lemma}
\label{quintlemma}
Let $w_i (x,t) $ be functions of space-time. Then
\begin{equation}
  \label{quint}
 \left| \int_0^\delta \int \prod_{i=1}^5 w_i (x,t) dx dt \right| \lesssim 
\left( \prod_{j=1}^3  {{\| w_j \|}_{X_{\frac{1}{4} , \half +}^\delta}} \right)
{{\| w_4 \|}_{X_{-\frac{3}{4} , \half +}^\delta}}
{{\| w_5 \|}_{X_{-\frac{3}{4} , \half +}^\delta}} .
\end{equation}
\end{lemma} %quintlemma

\begin{proof}
  The left-side of \eqref{quint} is estimated via H\"older's inequality
by
\begin{equation*}
  \left( \prod_{j=1}^3 {{\| w_j \|}_{L^4_x  L^\infty_{t \in [0,\delta]}}} 
\right)  {{\| w_4 \|}_{L^8_x  L^2_{t \in [0,\delta]}}}
{{\| w_5 \|}_{L^8_x  L^2_{t \in [0,\delta]}}}.
\end{equation*}
The first three factors are bounded using a maximal inequality  
from \cite{KPVCPAM},
\begin{equation}
  \label{fourmaximal}
  {{\| w \|}_{L^4_x  L^\infty_{t \in [0, \delta]}}} \lesssim {{\| w \|}_{
X^\delta_{\frac{1}{4}, \half+}}}.
\end{equation}
(Strictly speaking, \cite{KPVCPAM} contains an estimate for $S(t) \phi $
which implies \eqref{fourmaximal} by summing over cubic levels using
$b = \half+$, see \cite{B1} or \cite{GinibreNotes}, \cite{GinibreAsterisque}. 
A similar comment applies to \eqref{kato} below.)
The $w_4, ~w_5$ terms are controlled using the smoothing estimate
\begin{equation}
  \label{interpsmooth}
  {{\| w \|}_{L^8_x  L^2_t }} \leq 
{{\| w \|}_{X^\delta_{-\frac{3}{4}, \half+}}}
\end{equation}
which is an interpolant between the local-in-time energy estimate
\begin{equation}
  \label{energyest}
  {{\| w \|}_{L^2_x  L^2_t \in [0,\delta]}} \lesssim {{\| u \|}_{X^\delta_{0,
\half+}}}
\end{equation}
and the Kato smoothing estimate \cite{KPVCPAM}, 
valid for ${\widehat{w}}$ supported outside
$\{|\xi | < 1 \}$,
\begin{equation}
  \label{kato}
  {{\| w \|}_{L^\infty_x  L^2_{t \in [0, \delta]} }} 
\lesssim {{\| w \|}_{X^\delta_{-1, \half+}}}.
\end{equation}
In the remaining low frequency cases (e.g. when $\widehat{w_4}$ is 
supported inside $[-1,1]$) we have that
${{\| w_4 \|}_{L^\infty_x L^\infty_{t \in [0, \delta]}}} \leq  
{{\| w_4 \|}_{X^\delta_{0, \half+}}}$ and therefore may easily control
${{\| w_4 \|}_{L^8_x L^2_{t \in [0, \delta]}}}$ by 
$ {{\| w \|}_{X^\delta_{-\frac{3}{4}, \half+}}}$.
\end{proof}

Lemma \ref{quintlemma} is combined with the $M_5$ upper bound of 
Lemma \ref{MFivelemma} in the next result.

\begin{lemma}
  \label{incrementlemma}
Recall the definition \eqref{Ioperator} 
of the operator $I$. If the associated multiplier $m$ 
is of the form \eqref{particularm} with $s = -\frac{3}{4}+$ then
\begin{equation}
  \label{increment}
  \left| \int_0^\delta \Lambda_5 ( M_5; u_1, \dots, u_5 ) dt \right|
\lesssim N^{-\beta} \prod\limits_{j=1}^5 \Xd 0 {\half +} {Iu_j} ,
\end{equation}
with $\beta = 3 + \frac{3}{4} -.$
\end{lemma}

\begin{proof}
We may assume that the functions $\widehat{u_j}$ are nonnegative. By a
Littlewood-Paley decomposition, we restrict each $\widehat{u_j}$ to a 
frequency band $|\xi_j | \thicksim N_j$ (dyadic) and sum in the $N_j$
at the end of the argument.
The definition of the operator $I$ and \eqref{MFivebound} shows that it
suffices to prove
\begin{equation*}
\left|  \int_0^\delta \Lambda_5 \left( \frac{N_{45}  ~ m^2 (N_{*45} ) }{
(N+N_1 )(N+N_2) (N+N_3 ) (N+ N_{45} ) m(N_1) \dots m(N_5) } 
; u_1, \dots, u_5\right) dt \right| 
\end{equation*}
\begin{equation*}
\lesssim N^{-\beta} \prod_{j=1}^5 
N_j^{0-} {{\| u_j \|}_{X^\delta_{0, \half +}}}.
\end{equation*}
We cancel $\frac{N_{45}}{N+N_{45}} \leq 1$ and consider the
worst case when $m^2 ( N_{*45} ) =1 $ throughout. 

Note that $M_4$ vanishes when $|\xi_i | \ll N$ for $i = 1,2,3,4.$ Hence,
we are allowed to assume at least one, and hence two, of the $N_i \gtrsim N$.
Symmetry allows us to assume $N_1 \geq N_2 \geq N_3$ and $N_4 \geq N_5$.

The objective here is to show that 
\begin{equation*}
\left|  \int_0^\delta \Lam_5 \left(\prod_{i=1}^3 \frac{1}{(N+N_i ) m(N_i )}
\frac{1}{m(N_4 )} \frac{1}{m(N_5 )}; u_1, \dots, u_5 
\right) dt \right| \lesssim
N^{- \frac{15}{4} + } \prod_{j=1}^5 N_j^{0-} 
{{\| u_j \|}_{X^\delta_{0 , \half+}}} 
\end{equation*}

The form \eqref{particularm} of $m$ with $s = - \frac{3}{4}+$ implies
that $ \frac{1}{(N+N_i ) m(N_i )} \lesssim N^{-\frac{3}{4}+} 
\langle N_i \rangle^{-\frac{1}{4}-}.$ Therefore, we need to control
\begin{equation*}
  N^{-\frac{9}{4}+} \int_0^\delta \Lam_5 \left ( \langle N_1 
\rangle^{- \frac{1}{4} - }
\langle N_2 \rangle^{- \frac{1}{4} - } \langle N_3 
\rangle^{- \frac{1}{4} - } \frac{1}{m(N_4 )}
\frac{1}{m(N_5 )} \right) dt.
\end{equation*}
 
We break the analysis into three main cases: 
{\bf{Case 1.}} $N_{4} , N_5 \gtrsim N$,  
{\bf{Case 2.}} $N_4 \gtrsim N \gg N_5$, {\bf{Case 3.}} $N \gg N_4 \geq N_5$.

In Case 1, we have that $\frac{1}{m( N_4 )} \thicksim N^{-\frac{3}{4}+}
\langle N_4 \rangle^{\frac{3}{4}-}$ and $\frac{1}{m( N_5 )} 
\thicksim N^{-\frac{3}{4}+}
\langle N_5 \rangle^{\frac{3}{4}-}$. 
The desired prefactor $N^{-\frac{15}{4} +}$ then
appears and \eqref{quint} gives the result claimed.

In Case 2, $m(N_5) = 1$ and we must have $N_1 \geq N \geq N_5$ so
we multiply by ${{\left( \frac{N_1}{N} \right)}^{\frac{3}{4}} } 
{{\left( \frac{N_1}{N_5} \right)}^{\frac{1}{4}} } \geq 1$ and it suffices
to bound
\begin{equation}
\label{finished}
  N^{-\frac{3}{4}} N^{-\frac{12}{4} +} \left| \int_0^\delta \Lam_5 \left(
\langle N_1 \rangle^{\frac{3}{4}} \langle N_2  
\rangle^{-\frac{1}{4}-}  \langle N_3 \rangle^{-\frac{1}{4}-}
\langle N_4 \rangle^{\frac{3}{4}} \langle N_5 \rangle^{-\frac{1}{4}-} 
\right) dt \right|
\end{equation}
which may be done using \eqref{quint}.

For Case 3, we have $m(N_4) = m(N_5) =1.$ We are certain to have $N_1 \geq N_2
\geq N$ and can therefore multiply by
\begin{equation*}
  \frac{N_1}{N^{\frac{3}{4}} \langle N_4 \rangle^{\frac{1}{4}}}  
\frac{N_2}{N^{\frac{3}{4}} \langle N_5 \rangle^{\frac{1}{4}}} \geq 1
\end{equation*}
to again encounter \eqref{finished}.

\end{proof}

A glance back at \eqref{Efourinc} shows that for solutions of KdV, we
can now control the increment of the modified energy $E^4_I$.

\section{Global wellposedness of $KdV$ on $\R$}

The goal of this section is to construct the solution of the initial value
problem \eqref{KdV} on an arbitrary fixed time interval $[0,T]$. We 
first state a variant of the local wellposedness result of 
\cite{KPVBilin}. Next, we 
perform a rescaling under which the variant local result has an existence interval of size 1 and the initial data is small. This rescaling is possible
because the scaling invariant Sobolev index for $KdV$ is $- \frac{3}{2}$
which is much less than $-\frac{3}{4}$. Under the rescaling, we show that
\eqref{Efourinc} and \eqref{increment} allow us to iterate the local result
many times with an existence interval of size 1, thereby extending the 
local-in-time result to a global one. This will prove Theorem 1.

\subsection{A variant local wellposedness result}

The expression $\Lp 2 {Iu(t)}$, where ${\widehat{Iu(t)}}( \xi ) = m(\xi )
{\widehat{u(t)}}(\xi)$, and $m$ is of the form \eqref{particularm}, is
closely related to the $H^s ( \R )$ norm of $u$. Recall that the
definition of $m$ in \eqref{particularm} depends upon $s$. An adaptation
of the local well-posedness result in \cite{KPVBilin}, along the lines
of Lemma 5.2 in \cite{CKSTTDNLS1} and Section 12 in \cite{CKSTTGKdV}, 
establishes the following result.

\begin{proposition}
\label{variantlwp}
  If $s> - \frac{3}{4}$, the initial value problem \eqref{KdV} is locally
well-posed for data $\phi$ satisfying $I \phi \in L^2 (\R)$. Moreover,
the solution exists on a time interval $[0, \delta]$ with the lifetime
\begin{equation}
  \label{lifetime}
  \delta \thicksim {{\| I \phi \| }_{L^2 }^{-\alpha}},~\alpha >0,
\end{equation}
and the solution satisfies the estimate
\begin{equation}
  \label{spacetimebound}
  {{\| Iu \|}_{X^\delta_{0, \half+}}} \lesssim {{\| I \phi \|}_{L^2}}.
\end{equation}
\end{proposition}

\subsection{Rescaling}

Our goal is to construct the solution of \eqref{KdV} on an arbitrary
fixed time interval $[0, T]$. We rescale the solution by writing
$u_\lambda ( x, t) = \lam^{-2} u ( \frac{x}{\lam} , \frac{t}{\lam^3 } )$.
We acheive the goal if we construct $u_\lam $ on the time interval
$[0 , \lam^3 T ]$. A calculation shows that
\begin{equation*}
\Lp 2 { I \phi_\lam } \lesssim \lam^{- \frac{3}{2} - s } N^{-s} \Hsup s \phi.
\end{equation*}
The choice of the parameter $N= N(T) $ will be made later 
but we select $\lam$ now
by requiring
\begin{equation}
\label{chooselam}
\lam^{-\frac{3}{2} -s} N^{-s} \Hsup s \phi = \epsilon_0 < 1 \implies
\lam \thicksim N^{- \frac{2s}{3 + 2s}} .
\end{equation}
We drop the $\lam$ subscript on $u$ so that
\begin{equation}
\label{rescaled}
\Lp 2 {I \phi } = \epsilon_0 < 1
\end{equation}
and the task is to construct the solution of \eqref{KdV} on the time interval
$[0, \lam^3 T]$.

\begin{remark}
  The spatial domain for the initial value problem \eqref{KdV} is $\R$
which is invariant under the rescaling $x \longmapsto \frac{x}{\lam}.$
In contrast, the spatial domain $\T$ for the periodic initial value problem
for $KdV$ scales with $\lam$. The adaptation of our proof of global
well-posedness in the periodic context presented in Section 8 
requires us to identify the dependence of various estimates on the 
spatial period.
\end{remark}

\subsection{Almost conservation}
Recall the modified energy $E^2_I (0) = {{\| I \phi \|}_{L^2}^2}
= \Lam_2 ( m(\xi_1 ) m(\xi_2 )) (0).$ This subsection shows that the 
modified energy $E^2_I (t) $ of our rescaled local-in-time solution $u$ is
comparable to the modified energy $E^4_I (t)$. Next, as forecasted in
Section 5, we use \eqref{Efourinc} and the bound \eqref{increment} to show
$E^4_I (t)$ is almost conserved, implying almost conservation of
$E^2_I (t) = {{\| I u(t) \|}_{L^2}^2}$. Since the lifetime of the
local result \eqref{lifetime} is controlled by ${{\| I \phi \|}^{2}_{L^2}}$,
this conservation property permits us to iterate the local result with the 
same sized existence interval.

\begin{lemma}
  Let $I$ be defined with the multiplier $m$ of the form \eqref{particularm}
and $s = - \frac{3}{4}+.$ Then
\begin{equation}
\label{twobyfour}
| E^4_I (t) - E^2_I (t) | \lesssim {{\| Iu(t) \|}_{L^2}^3} + 
{{\| Iu(t) \|}_{L^2}^4}.
\end{equation}
\end{lemma}

\begin{remark}
  The estimate \eqref{twobyfour} is an a priori estimate for functions
of $x$ alone. The variable $t$ appears as a parameter.
\end{remark}

\begin{proof}
Since $E^4_I (t) = E^2_I (t) + \Lambda_3 ( \sigma_3 ) + 
\Lambda_4 (\sigma_4),$  it suffices to prove
\begin{equation}
  \label{threedsmall}
  |\Lambda_3 ( \sigma_3; u_1 , u_2 , u_3  ) | \lesssim \prod_{j=1}^3 
{{\| I u_j (t) \|}_{L^2}},
\end{equation}
\begin{equation}
  \label{fourdsmall}
  |\Lambda_4 ( \sigma_4 ; u_1, \dots, u_4 ) | \lesssim \prod_{j=1}^4 
{{\| I u_j (t) \|}_{L^2}}.
\end{equation}

We may again assume that the $\widehat{u_j}$ are nonnegative. 
By the definitions of $\sigma_3 $ \eqref{sigmathreedefined},
and $I$ \eqref{Ioperator}, and also \eqref{factthree} and \eqref{Mthreeviam},
\eqref{threedsmall}
follows if we show
\begin{equation}
\left| \Lambda_3 \left( \frac{ m^2 ( \xi_1 ) \xi_1 +  m^2 ( \xi_2 ) \xi_2
 + m^2 ( \xi_3 ) \xi_3}{ \xi_1 \xi_2 \xi_3 m( \xi_1 ) m( \xi_2 ) m( \xi_3 ) }
; u_1, u_2, u_3 \right) \right| \lesssim \prod_{j=1}^3 {{\| u_j \|}_2} .
\label{renormed}
\end{equation}
We make a Littlewood-Paley decomposition and
restrict attention
to the contribution arising from $|\xi_i | \thicksim N_i $ (dyadic),
and without loss assume $N_1 \geq N_2 \geq N_3$. In case $N_1 < \half N$,
then $m^2 ( \xi_i ) = 1,~i = 1,2,3 \implies \Lam_3 = 0.$ So, we
can assume $N_1 \thicksim N_2 \geq N_3 $. We consider separately the
cases: $N_3 \ll N,~N_3 \gtrsim N$.

{\bf{I.}} $N_3 \ll N$.  \newline
Since $\xi_1 + \xi_2 + \xi_3 =0 $ and $m^2$ controls itself (recall Lemma
\ref{onemvt}), we
may apply \eqref{mvt} to show $|m^2 ( \xi_1 ) \xi_1 + m^2 ( \xi_2 ) \xi_2
+ m^2 (\xi_3 ) \xi_3| \thicksim N_3$. Of course $m(N_3 ) =1$ in this case
so we need to bound $\Lambda_3 ( \frac{ N^s }{N_1^{1+s}} 
\frac{N^s}{N_1^{1+s}} )$. But this quantity is bounded by 
$\Lambda_3 (N_1^{-\frac{1}{6}} N_2^{-\frac{1}{6}} N_3^{-\frac{1}{6}} )$
(in fact with a decay in $N$) and we wish to prove
\begin{equation}
  \label{Holderdetails}
  \int\limits_{\xi_1 + \xi_2 + \xi_3 = 0,~|\xi_i | \thicksim N_i}
\prod_{i=1}^3  N_i^{-\frac{1}{6}} {\widehat{u_i}} (\xi_i ) \lesssim
\prod_{i=1}^3 {{\| u_i \|}_{L^2}}.
\end{equation}
Let $w_i ( x) $ be defined via
\begin{equation}
  \label{wdefined}
  {\widehat{w_i}} (\xi ) = N_i^{\frac{1}{6}} \widehat{u_i} (\xi )
\chi_{\{ |\xi | \thicksim N_i \}} ( \xi ).
\end{equation}
The left-side of \eqref{Holderdetails} may be rewritten
\begin{equation}
  \label{asconvolution}
  \int\limits_{\xi_3} {\widehat{w_3}} ( -\xi_3 ) \int\limits_{\xi_3 =
\xi_1 + \xi_2 } \widehat{w_1} (\xi_1 ) \widehat{w_2 } (\xi_2 )
= \langle \widehat{\overline{w_3} } , \widehat{u_1} \widehat{u_2} \rangle
\end{equation}
\begin{equation}
\label{asproduct}
= \langle \overline{w_3} , w_1 w_2 \rangle = \int w_3 w_1 w_2 dx.
\end{equation}
We may now apply H\"older in $L^3_{x} L^3_{x} L^3_{x}$ to bound
the left-side of \eqref{Holderdetails} by
\begin{equation*}
  {{\| v_3 \|}_{L^3_x}} {{\| w_1 \|}_{L^3_x }} {{\| w_2 \|}_{L^3_x}}.
\end{equation*}
Finally, the form of $w_i$ (and hence $v_3$) given in \eqref{wdefined}
allows us to conclude using Sobolev that
\begin{equation*}
  {{\| w_i \|}_{L^3_x}} \lesssim {{\| u_i \|}_{L^2_x}}.
\end{equation*}

\begin{remark}
  The argument reducing the left-side of \eqref{Holderdetails} to
\eqref{asproduct} by passing through the convolution representation
\eqref{asconvolution} will appear many times below. We will often compress
this discussion by referring to it as an ``$L^3_x L^3_x L^3_x$ H\"older
application''.
\end{remark}

{\bf{II.}} $N_3 \geq N.$ \newline
By definition of $m$, we have
\begin{equation*}
|m^2 ( \xi_1 ) \xi_1 + m^2 ( \xi_2 ) \xi_2 + m^2 ( \xi_3 ) \xi_3 |
\lesssim N^{-2s} ( N_1^{1+2s} + N_2^{1+2s} + N_3^{1+2s} ).
\end{equation*}
Suppose $s= - \frac{3}{4}+ < - \half$, then this expression is
\begin{equation*}
\thicksim N^{-2s} N_3^{1+2s}
\end{equation*}
Therefore, the multiplier in \eqref{renormed} is bounded by
\begin{equation*}
\frac{N_3^{1+2s}  N^{-2s} N^{3s} }{N_1^{1+s} N_2^{1+s} N_3^{1+s}}
\lesssim \frac{ N_3^s N^s }{ N_1^{1+s} N_2^{1+s}} \thicksim
N^s \frac{N_3^{s+ \frac{1}{6}} }{N_1^{1+s} N_2^{1+s} N_3^{\frac{1}{6}} }.
\end{equation*}
\begin{equation*}
\lesssim
N^{- \frac{3}{2} + \frac{1}{6}} N_1^{-\frac{1}{6}} N_2^{-\frac{1}{6}}
 N_3^{-\frac{1}{6}} 
\end{equation*}
and $L^3~ L^3 ~ L^3$ H\"older finishes off \eqref{renormed} and establishes
\eqref{threedsmall}.

We record here that the preceding calculations imply
\begin{equation}
\left| \frac{ m^2 ( \xi_1 ) \xi_1  + m^2 ( \xi_2 ) \xi_2  
+ m^2 ( \xi_3 ) \xi_3 }{
\xi_1 \xi_2 \xi_3 m( \xi_1 ) m( \xi_2 ) m( \xi_3 ) } \right|
\lesssim N^{-\frac{5}{4} +} N_1^{-\frac{1}{4}+} 
N_2^{-\frac{1}{4}+}  N_3^{-\frac{1}{4}+} .
\end{equation}

We turn our attention to proving \eqref{fourdsmall}. By \eqref{MFourEst} 
\eqref{factfour}, and the definition of $\sigma_4 $ \eqref{sigmafourdefined},
it suffices to control for $|\xi | \thicksim N_i$ (dyadic), with
$N_1 \geq N_2 \geq N_3 \geq N_4 \implies N_1 = N_2,$ that
\begin{equation}
\label{renormedfour}
\Lambda_4 \left( \frac{ 1 }{ (N+N_1 )(N+N_2 ) (N+ N_3 ) (N+ N_4 ) m( N_1 )
m( N_2 ) m( N_3) m( N_4) }; u_1, u_2, u_3 , u_4 \right) 
\end{equation}
\begin{equation*}
  \lesssim \prod_{j=1}^4 {{\| u_j \|}_{L^2}}.
\end{equation*}

The definition of $m$ shows the multiplier appearing in the left-side
of \eqref{renormedfour} is
\begin{equation*}
\lesssim \frac{N^{4s}}{ N_1^{1+s} N_2^{1+s} N_3^{1+s} N_4^{1+s} }
\end{equation*}
and for $s = - \frac{3}{4}+$,
\begin{equation}
\lesssim N^{4s} \frac{1}{N_3^{\half+} N_4^{\half+}} .
\end{equation}
With this upper bound on the multiplier, we bound the left-side of 
\eqref{renormedfour} in $L^2 ~ L^2 ~ L^\infty L^\infty $ via H\"older
and Sobolev to obtain the estimate 
\eqref{renormedfour} and therefore \eqref{fourdsmall}.
\end{proof}

Since our rescaled solution satisfies ${{\| I \phi \|}_{L^2}^2}
= \epsilon_0^2 < 1$, we are certain that
\begin{equation*}
  E^4_I (0) = E^2_I ( 0) + O ( \epsilon_0^3 ),
\end{equation*}
and, moreover, that
\begin{equation}
\label{comparable}
  E^4_I (t) = E^2_I (t) + O ( \epsilon_0^3 )  
\end{equation}
whenever ${{\| I u(t) \|}_{L^2}^2} = E^2_I (t) < 2 \epsilon_0 .$
Using the estimate \eqref{increment} in \eqref{Efourinc}, the rescaled solution
is seen to satisfy
\begin{equation}
  \label{Efourcontrol}
  E^4_I (t) \leq E^4_I (0) + C \epsilon_0^5 N^{-3 - \frac{3}{4}+} ~{\mbox{for all}}~t \in [0,1].
\end{equation}
Consequently, using \eqref{comparable}, we see that the rescaled
solution has 
\begin{equation*}
{{\| I u (1) \|}_{L^2}^2} = \epsilon_0^2 +  O (\epsilon_0^3)
+ C \epsilon_0^5 N^{-3 - \frac{3}{4}+} < 4 \epsilon_0^2 
\end{equation*}

\subsection{Iteration} 

We may now consider the initial value problem for KdV with initial data 
$u(1)$ and, in light of the preceding bound, the local result will advance
the solution to time $t=2$. We iterate this process $M$ times and, in place
of \eqref{Efourcontrol}, we have
\begin{equation*}
  E^4_I (t) \leq  E^4_I (0) + M C \epsilon_0^5 N^{- 3 - \frac{3}{4}+} ~{\mbox{for all}}~ t \in [0, M+1].
\end{equation*}
As long as $M N^{-3 - \frac{3}{4}+} \lesssim 1$, we will have the bound
  \begin{equation*}
    {{Iu(M)}_{L^2}^2} = \epsilon_0^2 + O (\epsilon^3 ) + M C \epsilon_0^5 
N^{-3 - \frac{3}{4}+} < 4 \epsilon_0^2,
  \end{equation*}
and the lifetime of the local results remain uniformly of size 1. We take
$M \thicksim N^{3 + \frac{3}{4}-}.$
This process extends the local solution to the time interval 
$[0, N^{3 + \frac{3}{4}-} ]$. We choose $N = N(T)$ so that 
\begin{equation*}
  N^{3 + \frac{3}{4} - } > \lam^3 T \thicksim N^{- \frac{6s}{3 + 2s} } T,
\end{equation*}
which may certainly be done for $s > - \frac{3}{4}.$ This completes the
proof of global well-posedness for $KdV$ in $H^s (\R ), ~ s> - \frac{3}{4}$.

We make two observations regarding the rescalings of our global-in-time 
KdV solution:
\begin{equation}
\label{polyone}
\sup_{t \in [0,T]} \Hsup s {u(t) } \thicksim \lam^{\frac{3}{2} + s}
\sup_{t \in [0 ,\lam^3 T]} \Hsup s {u_\lam (t) } \leq \lam^{\frac{3}{2} + s}
\sup_{t \in [0, \lam^3 T]} {{ \| I u_\lam (t) \|}_{L^2}},
\end{equation}
\begin{equation}
\label{polytwo}
{{\| I \phi_\lam \|}_{L^2}} \lesssim N^{-s } \Hsup s {\phi_\lam } \thicksim
N^{-s} \lam^{- \frac{3}{2} -s } \Hsup s \phi .
\end{equation}

The almost conservation law and local well-posedness iteration argument
presented above implies that {\it{provided}} $N$ {\it{and}} $\lam$ 
{\it{are selected correctly}}
\begin{equation}
  \label{polythree}
  \sup_{t \in [0, \lam^3 T]} {{\| I u_\lam (t) \|}_{L^2}} \lesssim
  \sup_{t \in [0, \lam^3 T]} {{\| I \phi_\lam \|}_{L^2}}.
\end{equation}
The estimate \eqref{polythree} forms a bridge between \eqref{polyone}
and \eqref{polytwo} which implies
\begin{equation}
  \label{polyfour}
  \sup_{t \in [0,T]} \Hsup s {u(t) } \lesssim N^{-s} \Hsup s {\phi }.
\end{equation}
In fact, the selection of $N$ is polynomial in the parameter $T$ so
\eqref{polyfour} gives a polynomial-in-time upper-bound on $\Hsup s {u(t)}$.

{\bf{The choice of $\lam$}} \newline
The parameter $\lam$ was chosen above so that 
\begin{equation}
\label{polyfive}
{{\| I \phi_\lam \|}_{L^2}}
\thicksim \epsilon_0 \ll 1.
\end{equation}
 Since, from \eqref{polytwo}, ${{\| I \phi_\lam \|}_{L^2}} \lesssim N^{-s} \lam^{-\frac{3}{2} -s } \Hsup s \phi $, we see that \eqref{polyfive} holds
provided we choose
\begin{equation}
  \label{polysix}
  \lam = \lam (N, \epsilon_0 , \Hsup s \phi ) \thicksim
{{\left( \frac{ \Hsup s \phi }{\epsilon_0 } \right)}^{\frac{2}{3+2s}}}
N^{-\frac{2s}{3+2s}} .
\end{equation}

{\bf{The choice of $N$}} \newline
The parameter $N$ is chosen so that
\begin{equation}
  \label{polyseven}
  N^\beta > \lam^3 T \thicksim c_{\Hsup s \phi , \epsilon_0 } 
N^{-\frac{6s}{3+2s}} T ,
\end{equation}
where $\beta$ is the exponent appearing in \eqref{increment} (in the $\R$-case
just presented, $\beta = 3 + \frac{3}{4} -$). This unravels to give
a sufficient choice of $N$:
\begin{equation}
  \label{polyeight}
  N \thicksim  c_{\Hsup s \phi , \epsilon_0 } T^{ \frac{3 + 2s}{\beta (3 + 2s)
+ 6s }} \thicksim   c_{\Hsup s \phi , \epsilon_0 } T^{\gamma(s)}.
\end{equation}
In the range $-\frac{3}{2} < s$, the numerator of the exponent on $T$ is
positive. The denominator is positive provided $\beta > - \frac{6}{3+2s}$.
For $s = -\frac{3}{4},~- \frac{6}{3+2s} = 3$ so we require better than
third order decay with $N$ in the local-in-time increment \eqref{increment}.
With $s = -\frac{3}{4}+, ~\beta = 3 + \frac{3}{4} -$, calculating
$\gamma (s)$ and inserting the resulting expression for $N$ in terms of
$T$ into \eqref{polyfour} reveals that, for our global-in-time solutions of
\eqref{KdV}, we have
\begin{equation}
  \label{polyboundonR}
  {{\| u(t) \|}_{H^{-\frac{3}{4}+} (\R )}} \lesssim t^{1+} 
{{\| \phi \|}_{H^{-\frac{3}{4}+}}} .
\end{equation}

\begin{remark}
Observe that the polynomial exponent $1+$ in \eqref{polyboundonR} does
not explode as we approach the critical regularity value $- \frac{3}{4}$.
This is due to the fact that \eqref{increment} gave us much more decay
than required for iterating the local result. In principle, the
decay rate in \eqref{increment} could be improved by going further
along the sequence $\{ E^n_I \}$ of modified energies. If local
well-posedness of KdV is proven in $H^{-\frac{3}{4} } (\R )$, the
bounds obtained here should give global well-posedness in $H^{-\frac{3}{4}}
(\R )$.
\end{remark}

\section{Local well-posedness of KdV on $\T$}

This section revisits the local-in-time theory for periodic KdV
developed by Kenig, Ponce and Vega \cite{KPVBilin} and Bourgain \cite{B1}.
Our presentation provides details left unexposed in \cite{KPVBilin} and 
\cite{B1} and quantifies the dependence of various implied constants on the 
length of the spatial period. This quantification is necessary for the 
adaptation of the rescaling argument used in Section 6 to the
periodic setting.

\subsection{The $\lambda$-periodic initial value problem for KdV}

We consider the {\it{$\lambda$-periodic initial value problem for KdV}}:
\begin{equation}
  \label{lamkdv}
  \left\{
   \begin{matrix}
    \partial_t u + \partial_x^3 u + \half \partial_x u^2 =0,& x \in [0, \lam] 
         \\
     u(x, 0) = \phi (x).
   \end{matrix}
\right.
\end{equation}
We first want to build a representation formula for the solution of the 
linearization of \eqref{lamkdv} about the zero solution. So, we wish
to solve the {\it{linear homogeneous $\lambda$-periodic initial value 
problem}} 
\begin{equation}
  \label{lamhomog}
  \left\{
   \begin{matrix}
    \partial_t w + \partial_x^3 w   =0,& x \in [0, \lam] 
         \\
     w(x, 0) = \phi (x).
   \end{matrix}
\right.
\end{equation}
Define $(dk)_\lam$ to be normalized counting measure on $\Z / \lam $:
\begin{equation}
\label{dklam}
\int a(x) (dk )_\lam = \frac{1}{\lam} \sum_{k \in \Z / \lam } a(k).
\end{equation}
Define the {\it{Fourier transform}} of a function $f$ defined on $[0, \lam]$
by 
\begin{equation}
  \label{FourierTransform}
  {\widehat{f}} ( k) = \int_0^\lam e^{-2 \pi i k x } f(x) dx
\end{equation}
and we have the {\it{Fourier inversion formula}}
\begin{equation}
\label{FourierInversion}
f(x) = \int e^{2 \pi i k x } {\widehat{f}} (k) (dk)_\lam .
\end{equation}
The usual properties of the Fourier transform hold:
\begin{equation}
  \label{Plancherel}
  {{\| f \|}_{L^2 ( [0,\lam ] )}} = 
{{\| \widehat{f} \|}_{L^2 ( ( dk )_\lam )}} ~~~~{\mbox{(Plancherel)}},
\end{equation}
\begin{equation}
  \label{Parseval}
  \int_0^\lam  f(x) {\overline{g(x)}} dx = \int {\widehat{f}}( k ) 
{\overline{\widehat{g}} (k) } (dk)_\lam ~~~~{\mbox{(Parseval)}},
\end{equation}
\begin{equation}
  \label{Convolution}
  {\widehat{fg}} (k) = {\widehat{f}} *_\lam {\widehat{g}} (k) =
\int {\widehat{f}} ( k - k_1 ) {\widehat{g}} (k_1 ) (d k_1 )_\lam 
~~~~{\mbox{(Convolution)}},
\end{equation}
and so on.
If we apply $\partial_x^m, ~ m \in \N$ to \eqref{FourierInversion} we obtain
\begin{equation*}
  \partial_x^m f(x) = \int e^{2 \pi i k x } {{(2 \pi i k )}^m} {\widehat{f}}
(k ) (dk )_\lam.
\end{equation*}
This, together with \eqref{Plancherel}, motivates us to define the 
{\it{Sobolev space}} $H^s ( 0 , \lam )$ with the norm
\begin{equation}
  \label{Hs}
  {{\| f \|}_{H^s ( 0 , \lam )}} = {{\| {\widehat{f}} ( k ) \langle k 
\rangle^s \|}_{L^2 ( 
(dk)_\lam )}}.
\end{equation}
We will often denote this space by $H^s$ for simplicity. 
Note that there are about $\lam$ {\it{low frequencies}}
in the range $|k| \lesssim 1$ where the $H^s$ norm consists of the $L^2$
norm.

The Fourier inversion formula \eqref{FourierInversion} allows us to write
down the solution of \eqref{lamhomog}:
\begin{equation}
  \label{Slam}
  w(x,t) = S_\lam (t) \phi (x) = \int e^{2 \pi i k x} e^{- {{(2 \pi i k)}^3} t}
{\widehat{\phi}} (k) (dk )_\lam. 
\end{equation}

For a function $v = v(x,t)$ which is $\lam$-periodic with respect to the 
$x$ variable and with the time variable $t \in \R$, we define the
{\it{space-time Fourier transform}} ${\widehat{v}} = {\widehat{v}}(k , \tau )$
for $k \in \Z / \lam$ and $\tau \in \R$ by
\begin{equation}
  \label{spacetimeFT}
  {\widehat{v}} (k, \tau ) = \int \int_0^\lam e^{-2 \pi i k x }
e^{- 2 \pi i \tau t } v(x,t) dx dt.
\end{equation}
This transform is inverted by
\begin{equation}
  \label{spacetimeinvFT}
  v(x,t) = \int \int e^{2 \pi i k x } e^{2 \pi i \tau t } {\widehat{v}}
( k , \tau ) (dk )_\lam d\tau .
\end{equation}
The expression \eqref{Slam} may be rewritten as a space-time inverse
Fourier transform,
\begin{equation}
  \label{Slamcurve}
  S_\lam (t) \phi (x) = \int \int e^{2 \pi i k x} e^{2 \pi i \tau t}
\delta ({  \tau - 4 \pi^2 k^3  }) {\widehat{\phi}} 
(k) (dk)_\lam d\tau 
\end{equation}
where $\delta ({ \eta  } )$ represents a 1-dimensional Dirac mass
at $\eta = 0$. This recasting shows that $S_\lam ( \cdot ) \phi $
has its space-time Fourier Transform supported precisely on the
cubic $\tau = 4 \pi^2 k^3 $ in $\Z / \lam \times \R$.

We next find a representation for the solution of the {\it{linear
inhomogeneous $\lambda$-periodic initial value problem}}
\begin{equation}
  \label{laminhomog}
  \left\{
\begin{matrix}
\partial_t v + \partial_x^3 v = f, & x \in [0, \lam ] \\
v(x, 0 ) = 0,
\end{matrix}
\right.
\end{equation}
with $f = f(x,t)$ a given time dependent $\lam$-periodic (in $x$) function.
By Duhamel's principle,
\begin{equation}
  \label{lamDuh}
  v(x,t) = \int_0^t S_\lam (t - t') f(x, t') dt'.
\end{equation}
We represent $f(x,t')$ using \eqref{spacetimeinvFT}, apply \eqref{Slam}
and rearrange integrations to find
\begin{equation*}
  v(x,t) = \int \int e^{2 \pi i k x } e^{2 \pi i (4 \pi^2 k^3 t )}
\int_0^t e^{2 \pi i (\tau - 4 \pi^2 k^3 ) t' } dt' {\widehat{f}}(k , \tau)
(dk)_\lam d\tau .
\end{equation*}
Performing the $t'$-integration, we find
\begin{equation}
  \label{afterintegration}
  v(x,t) = \int \int e^{2 \pi k x } e^{2 \pi i (4 \pi^2 k^3 t )}
\frac{e^{2 \pi i (\tau - 4 \pi^2 k^3) t} - 1}{2 \pi i (\tau - 4 \pi^2 
k^3 )}  {\widehat{f}}(k , \tau)
(dk)_\lam d\tau .
\end{equation}

The $\lam$-periodic initial value problem for KdV \eqref{lamkdv}
is equivalent to the integral equation
\begin{equation}
  \label{lamintegraleq}
  u(t) = S_\lam (t) \phi - \int_0^t S_\lam (t -t') ( \half \partial_x
u^2 ( t')) dt'.
\end{equation}

\begin{remark}
\label{meanzeroremark}
The spatial mean $\int_{\T} u(x,t) dx$ is conserved during the evolution
\eqref{lamkdv}. We may assume that the initial data $\phi$ satisfies a 
mean-zero assumption $\int_{\T} \phi (x) dx$ since otherwise we can replace
the dependent variable $u$ by $v = u - \int_{\T} \phi $ at the expense
of a harmless linear first order term. This observation was used by Bourgain
in \cite{B1}. The mean-zero assumption is crucial for some of the analysis
that follows.
\end{remark}

\subsection{Spaces of functions of space-time}

The integral equation \eqref{lamintegraleq} will be solved using the
contraction principle in spaces introduced in this subsection. We
also introduce some other spaces of functions of space-time which
will be useful in our analysis of \eqref{lamintegraleq}.

We define the $\Xsb$ spaces for $\lam$-periodic KdV via the norm
\begin{equation}
  \label{lamXsb}
  {{\| u \|}_{\Xsb ( [0, \lam ] \times \R)}} =
{{\| \langle k \rangle^s \langle \tau - 4 \pi^2 k^3 \rangle^b \widehat{u} (k , \tau ) \|}_{L^2 ( (dk)_\lam 
d \tau )}}.
\end{equation}
(We will suppress reference to the spatial period $\lam$ in the notation
for the space-time function spaces $X_{s,b}$ and the related spaces below.)
These spaces were first used to systematically study nonlinear
dispersive wave problems by Bourgain \cite{B1}. Klainerman and Machedon
\cite{KlainMach94} 
used similar ideas in their study of the nonlinear wave equation.
The spaces appeared earlier in a different setting in the work \cite{Beals}
of M. Beals. 

The study of periodic KdV in \cite{KPVBilin}, \cite{B1} has been based
around iteration in the spaces $X_{s, \half}$. This space barely fails to
control the $L^\infty_t H^s_x $ norm. To ensure continuity of the time
flow of the solution we construct, we introduce the slightly smaller space
$Y^s$ defined via the norm
\begin{equation}
  \label{Ys}
  {{\| u \|}_{Y^s}} = \X s \half u  + {{\|\langle k \rangle^s \widehat{u} (k , \tau) \|}_{
L^2 ( (dk)_\lam ) L^1 (d \tau )}}.
\end{equation}
If $u \in Y^s$ then $u \in L^\infty_t H^s_x $.
We will construct the solution of \eqref{lamintegraleq} by proving a 
contraction estimate in the space $Y^s$. The mapping properties of 
\eqref{lamDuh} motivate the introduction of the companion spaces $Z^s$ defined
via the norm
\begin{equation}
  \label{Zs}
  {{\| u \|}_{Z^s}} = \X s {-\half} u   + {{\left\| 
\frac{\langle k \rangle^s \widehat{u} (k , \tau)}{\langle \tau - 4 
\pi^2 k^3\rangle} \right\|}_{L^2 
((dk)_\lam ) L^1 (d\tau) }}.
\end{equation}

Let $\eta \in C_0^\infty (\R)$ be a nice bump function supported on 
$[-2,2]$ with $\eta =1$ on $[-1,1]$. It is easy to see that multiplication
by $\eta (t)$ is a bounded operation on the spaces $Y^s$, $Z^s$, and
$\Xsb$.

\subsection{Linear estimates}

\begin{lemma}
  \label{lamhomogest}
\begin{equation}
\label{lamhomogestimate}
{{\| \eta(t) S_\lam (t) \phi \|}_{Y^s}} \lesssim {{\| \phi \|}_{H^s}}.
\end{equation}
\end{lemma}

The proof follows easily from the fact that
\begin{equation}
  \label{tauconv}
  {\widehat{\eta  S_\lam  ( \phi) }} ( k, \tau ) = \widehat{\phi}
(k) {\widehat{\eta}} ( \tau - 4 \pi^2 k^3 ).
\end{equation}

\begin{lemma}
  \label{laminhomogest}
\begin{equation}
\label{laminhomogestimate}
{{\left\| \eta(t) \int_0^t S_\lam (t - t') F(t') dt' \right\|}_{Y^s}}
\lesssim {{\| F \|}_{Z^s}}.
\end{equation}
\end{lemma}

\begin{proof}
By applying a
smooth cutoff, we may assume that $F$ is supported on $\T \times [-3, 3]$.
Let $a(t) = {\mbox{sgn}} (t) {\tilde{\eta}} ( t )$, where ${\tilde{\eta}}$
is a smooth bump function supported on $[-10, 10]$ which equals 1 on
$[-5,5]$. The identity
\begin{equation*}
\chi_{[0,t]} ( t') = \half ( a( t') - a ( t-t') ),
\end{equation*}
which is valid for $t \in [-2,2]$ and $t' \in [-3,3]$, allows us to
rewrite $\eta (t) \int_0^t S ( t-t') F(t') dt'$ as a linear combination
of
\begin{equation}
\label{firstpart}
\eta(t) S(t) \int_{\R} a(t') S( - t') F(t') dt' 
\end{equation}
and 
\begin{equation}
\label{secondpart}
\eta(t) \int_{\R} a( t - t') S(t-t') F(t') dt'.
\end{equation}
Consider the contribution \eqref{firstpart}. By \eqref{lamhomogestimate},
it suffices to show that
\begin{equation*}
\Hsup s {\int a(t') S(-t') F(t') dt' }  \lesssim {{\| F \|}_{Z^s}} .
\end{equation*}
Since the Fourier transform of $\int a(t') S( - t') F(t') dt'$
evaluated at $\xi $ is given by $\int \widehat{a} ( \tau - \xi^3 )
\widehat{F} ( \xi , \tau ) d\tau $ and one can easily verify that
$|\widehat{a} ( \tau ) | = O ( \langle \tau \rangle^{-1} )$, 
the claimed estimate
follows using the definition \eqref{Zs}.

For \eqref{secondpart}, we discard the cutoff $\eta (t) $ and note
that the space-time Fourier transform of $\int a( t-t')  S(t-t' )
F (t') dt'$ evaluated at $(\xi , \tau )$ is equal to $\widehat{a}
( \tau - 4 \pi^2 \xi^3 ) \widehat{F} ( \xi , \tau )$. The claimed estimate
then follows from the definitions \eqref{Zs}, \eqref{Ys} and the
decay estimate for $\widehat{a}$ used above.
\end{proof}

\begin{lemma}
  \label{lamoncurveStrichartz}
Let $\phi$ be a $\lam$-periodic function whose Fourier Transform is supported
on $\{ k : |k| \thicksim N\}$. Then
\begin{equation}
\label{LFouroncurve}
{{\| \eta (t) S_\lam (t) \phi 
\|}_{L^4_{x,t}}} \lesssim C(N, \lam) {{\| \phi \|}_{L^2_x}},
\end{equation}
where
\begin{equation}
  \label{cofNlam}
  C(N, \lam ) = \left\{ \begin{matrix} 1 & {\mbox{if}}~ N \leq 1, \\
  \left( \frac{1}{\sqrt{N}} + \frac{1}{\lam} \right)^{\frac{1}{4}} 
& {\mbox{if}}~ N \geq 1.
\end{matrix}
\right.
\end{equation}
\end{lemma}
 
\begin{remark}
In the limit $\lam  \rightarrow \infty$, \eqref{LFouroncurve} yields
the Strichartz estimate on the line (at least when $N \geq 1$),
\begin{equation}
  \label{StrichEighth}
  {{\| D_x^{\frac{1}{8}} e^{-t 4 \pi^2 
\partial_x^3} \phi \|}_{L^4_{x \in \R ,t}}}
\lesssim {{\| \phi \|}_{L^2_{x \in \R}}}.
\end{equation}
\end{remark}

\begin{proof}
It suffices to show that
\begin{equation}
  \label{bilinearize}
  {{\| (\eta(t) S_\lam (t) \phi_1 ) (\eta(t) S_\lam (t) \phi_2) \|}_{
L^2_{x,t}}^2} \leq C^2 (N, \lam ) {{\| \phi_1 \|}_{L^2_x}} {{\| \phi_2 
\|}_{L^2_x}} 
\end{equation}
for functions $\phi_1, ~\phi_2$ satisfying the hypotheses. Properties of
the Fourier Tranform allow us to reexpress the left-side as
\begin{equation*}
  {{\left\| \int\limits_{k=k_1 +k_2 , \tau = \tau_1 + \tau_2 }
{\widehat{\phi_1 }} ( k_1 ) {\widehat{\phi_2}} ( k_2 ) \psi( \tau_1
- 4\pi^2 k_1^3 ) \psi (\tau_2 - 4\pi^2 k_2^3 ) (dk_1)_\lam d\tau_1 \right\|}_{
L^2 (d\tau ~ (dk)_\lam )}^2},
\end{equation*}
where $\psi = {\widehat{\eta}}$, may be take to be a positive even
Schwarz function. We evaluate the $\tau_1$-integration by writing
\begin{equation*}
  \int \psi( \tau_1 - 4 \pi^2 k_1^3 ) \psi ( \tau - \tau_1 - 4 \pi^2 
k_2^3 ) d\tau_1
= {\tilde{\psi}} ( \tau - 4 \pi^2 k_1^3 - 4 \pi^2 k_2^3 ),
\end{equation*}
with ${\tilde{\psi}}$ also rapidly decreasing. Inserting this into the
reexpressed left-side and applying Cauchy-Schwarz leads to the upper
bound
\begin{equation*}
  {{\left\| {{\left( \int {\tilde{\psi}}^2 ( \tau - 4 \pi^2 k_1^3 - 
4 \pi^2 k_2^3 )
(dk_1 )_\lam \right)}^\half} {{\left( \int  
{\tilde{\psi}}^2 ( \tau - 4 \pi^2 k_1^3 - 4 \pi^2 k_2^3 ) 
|\widehat{\phi_1} (k_1 )|^2 
 |\widehat{\phi_2} (k_2 )|^2  (dk_1 )_\lam \right)}^\half} \right\|}_{L^2 
(d\tau (dk)_\lam )}}.
\end{equation*}
The first integral may be pulled out of the $L^2$ norm and the $\tilde{\psi^2}$
term in the second integral is used to integrate in $\tau$ to give
\begin{equation*}
  \lesssim {{\| \int {\tilde{\psi}}^2 ( \tau - 4 \pi^2 k_1^3 - 4 \pi^2 
k_2^3 ) (dk_1 )_\lam 
\|}_{
L^\infty_{k , \tau}}} {{\| \phi_1 \|}_{L^2_x}} {{\| \phi_2 \|}_{L^2_x}}.
\end{equation*}
Matters are thus reduced to quantifying the $L^\infty$ norm above. Let 
$M$ denote  ${{\| \int {\tilde{\psi}}^2 ( \tau - 4 \pi^2 k_1^3 - 4 \pi^2 
k_2^3 ) (dk_1 )_\lam 
\|}_{
L^\infty_{k , \tau}}}$. We estimate $M$ by counting
\begin{equation*}
  M \lesssim  
\frac{1}{\lam} | \{ k_1 \in \Z / \lam : |k_1 | \thicksim N ; 
|k-k_1 | \thicksim N; k^3 - 3k k_1 (k-k_1 ) = \tau + O(1) \} |.
\end{equation*}
In case $N \leq 1$, the cardinality of the set is $O(\lam )$ so
$C(N,\lam )  \lesssim 1$ for $N \leq 1$. Assume now that $N > 1$
and rename $k_1 = x$. The task is to estimate
\begin{equation*}
  |\{ x \in \Z /\lam : |x|, |k -x | \thicksim N; 3k (x-\frac{k}{2})^2 - 
\frac{k^2}{4} = \tau - k^3 +O(1) \}|.
\end{equation*}
This set is largest when the parabola is the flattest, so when $x \thicksim
\frac{k}{2}$. We find that
\begin{equation*}
  M \lesssim \frac{1}{\lam } \left( \frac{1}{\sqrt{k}} \lam + 1 \right)
=  \left( \frac{1}{\sqrt{k}}  + \frac{1}{\lam} \right),
\end{equation*}
which completes the proof.
\end{proof}

\begin{lemma} If $v=v(x,t) $ is a $\lam$-periodic function of $x$
and the spatial Fourier tranform of $v$ is supported on $\{ k: |k| 
\thicksim N\}$ then
\begin{equation}
  \label{lamBouStrichartz}
  {{\|  \eta(t) v \|}_{L^4_{x,t}}} 
\lesssim C(N, \lam ) \X 0 {\half+} v .
\end{equation}
where $C(N, \lam)$ is as it appears in \eqref{cofNlam}.
\end{lemma}
This follows easily by stacking up cubic level sets on 
which \eqref{LFouroncurve} holds.

\begin{lemma}
If $v=v(x,t)$ is a $\lam$-periodic function of $x$ then
\begin{equation}
  \label{rescaledBourgain}
  {{\|\eta(t)  v \|}_{L^4_{x,t}}} \lesssim  \X 0 {\frac{1}{3}} v .
\end{equation}
\end{lemma}

The estimate \eqref{rescaledBourgain} is a rescaling of 
the $\lam =1$ case proven by Bourgain \cite{B1}.

\begin{remark}
We can interpolate between \eqref{lamBouStrichartz} and 
\eqref{rescaledBourgain} to obtain
\begin{equation}
  \label{atahalf}
  {{\| \eta (t) v \|}_{L^4_{x,t}}} \lesssim \{C(N, \lam)\}^{1-} \X 0 \half v .
\end{equation}
\end{remark}

\subsection{Bilinear estimate}

\begin{proposition}
\label{lambilinear}
If $u$ and $v$ are $\lam$-periodic functions of $x$, also depending upon $t$
having zero $x$-mean for all $t$, then
\begin{equation}
{{\| \eta (t) \partial_x ( u v) \| }_{Z^{-\half}}} \lesssim \lam^{0+} 
{{\| u \|}_{X_{-\half, \half}}} {{\| v \|}_{X_{-\half, \half}}}.
\label{bilinearestimate}
\end{equation}
\end{proposition}

Note that \eqref{bilinearestimate} implies ${{\| \eta (t) \partial_x ( u v) \| }_{Z^{-\half}}} \lesssim \lam^{0+} 
{{\| u \|}_{Y^{-\half}}} {{\| v \|}_{Y^{-\half}}}.$ We will relax the notation
by dispensing with various constants involving $\pi$ with the recognition that
some of the formulas which follow may require adjusting the constants to be
strictly correct.

\begin{proof} The norms involved allow us to assume that $\widehat{u}$
and $\widehat{v}$ are nonnegative.
There are two contributions to the $Z^{-\half}$ norm we must control.
We begin with the $X_{-\half, -\half}$ contribution.  
Duality and an integration by parts
shows that it suffices to prove
\begin{equation*}
  \int \int u_1 (x,t) u_2 (x,t) w_x (x,t) \eta(t) dx dt
\lesssim \lam^{0+} \X {-\half} {\half} {u_1}  \X {-\half} {\half} {u_2}
\X {\half } {\half } {w} .
\end{equation*}
Writing $u_3 = w_x$ shows that it suffices to prove
\begin{equation}
  \label{symreduced}
  \int 
\int u_1 (x,t) u_2 (x,t) u_3 (x,t) \eta(t) dx dt \lesssim \lam^{0+}
\X {-\half} {\half} {u_1} \X {-\half} {\half} {u_2} \X {-\half} {\half} {u_3} 
\end{equation}
for all $u_1, ~u_2, ~u_3$ having zero $x$-mean. The left-side may be 
rewritten
\begin{equation}
  \label{rewrite}
  \int\limits_{k_1 +k_2 +k_3 =0, } \int\limits_{\tau_1 + \tau_2 + \tau_3 =0} 
{\widehat{u_1}} ( k_1, \tau_1 ) {\widehat{u_2}} ( k_2, \tau_2 ) 
{\widehat{{u_3}}} ( k_3, \tau_3 ) (dk_1)_\lam (dk_3)_\lam
d\tau_1 d\tau_2 d\tau_3 .
\end{equation}
Note that the mean zero conditions allow us to assume $k_i \neq 0$.

{\bf{Case 1.}} $|k_1 |, |k_2 | , |k_3 | \gtrsim 1 $. \newline
The identity
\begin{equation*}
  3 k_1 k_2 k_3 = k_1^3 + k_2^3 +k_3^3
\end{equation*}
and the Case 1 defining conditions imply
\begin{equation*}
  1 \lesssim \sum\limits_{j=1}^3 
\frac{ |\tau_j - k_j^3 | }{\langle k_1 \rangle  \langle k_2 \rangle \langle 
k_3 \rangle }.
\end{equation*}
This of course implies
\begin{equation}
  \label{powerfuldenoms}
  1 \lesssim \sum\limits_{j=1}^3 
\frac{ \langle \tau_j - k_j^3 \rangle^{\half} }{\langle k_1 \rangle^{\half} \langle k_2 \rangle^{\half} 
\langle k_3 \rangle^{ \half} }.
\end{equation}
Inserting \eqref{powerfuldenoms} into \eqref{rewrite} and using symmetry
reduces matters to showing that
\begin{equation*}
   \int\limits_{k_1 +k_2 +k_3 =0 } \int \widehat{\eta} ( \tau_1 +\tau_2 
+\tau_3) {\frac{\widehat{u_1} (k_1 , \tau_1)}{{{\langle k_1 \rangle}^\half}}}
 {\frac{\widehat{u_2} (k_2 , \tau_2 )}{{{\langle k_2 \rangle}^\half}}}
   {\frac{\widehat{u_3} (k_3 , \tau_3 ) 
\langle \tau_3 - k_3^3 \rangle^{\half} }{{{\langle k_3 \rangle}^\half}}}
\end{equation*}
\begin{equation*}
  \lesssim \lam^{0+}  \X {-\half} {\half} {u_1}    \X {-\half} {\half} {u_2}
 \X {-\half} {\half} {u_3}
\end{equation*}
After some natural substitutions and undoing the Fourier Transform, 
we see that it suffices to show that
\begin{equation}
  \label{CaseOneReduction}
  \left| \int \int \eta(t) v_1 v_2 v_3 dx dt \right|
\lesssim \lam^{0+}  \X 0 \half {v_1}  \X 0 {\half} {v_2}  \X 0 0 {v_3} .
\end{equation}
We split the left-side using H\"older into $L^4_{x,t} L^4_{x, t} L^2_{x,t}$
and then apply \eqref{rescaledBourgain} to finish off this case. (In fact,
we control the left side of \eqref{CaseOneReduction} with $\lam^{0+}
\X 0 {\frac{1}{3}} {v_1}  \X 0 {\frac{1}{3}} {v_2}  \X 0 0 {v_3} .$)

{\bf{Case 2.}} $|k_1|, |k_2|, |k_3| \lesssim 1.$ \newline
Derivatives are cheap in this frequency setting.
We use H\"older to estimate \eqref{rewrite} in $L^4 L^4 L^2$,
then apply \eqref{rescaledBourgain} to control the $L^4$ norms. Finally, we
use the Case 2. defining conditions and Sobolev 
to move $X_{0,\frac{1}{3}}$ to $X_{-\half, \half}$ on two factors and
${X_{0,0}}$ to $X_{-\half, \half}$ on the remaining factor. (Again, we
have encountered $X_{0, \frac{1}{3}}$ on two
factors.)

Since $k_1 + k_2 + k_3 = 0$, the only remaining case to consider
is when one of the frequencies is small and the other two are big.
Symmetry permits us to focus on \newline
{\bf{Case 3.}} $0 < |k_3| \lesssim 1 \lesssim |k_1|, |k_2|.$ \newline

The analog of \eqref{powerfuldenoms} in this case is
\begin{equation}
  \label{lesspowerful}
  1 \lesssim |k_3|^{-\half} \sum\limits_{j=1}^3 
\frac{\langle \tau_j - k_j^3 \rangle^\half}{\langle k_1 \rangle^\half \langle 
k_2 \rangle^\half \langle k_3 \rangle^\half } .
\end{equation}
Since we are in the $\lam$-periodic setting and our functions have zero
$x$-mean, we have $|k_3 | \gtrsim \frac{1}{\lam}.$ Therefore, we may
replace \eqref{lesspowerful} by the symmetric expression
\begin{equation}
  \label{symlesspowerful}
  1 \lesssim \lam^\half \sum\limits_{j=1}^3 
\frac{\langle \tau_j - k_j^3 \rangle^\half}{\langle k_1 \rangle^\half \langle k_2 \rangle^\half \langle k_3 \rangle^\half } .
\end{equation}
Repeating certain arguments in Case 1 reduces matters to showing that
\begin{equation}
\label{CaseThreeReduction}
\lam^\half   \left| \int \int \eta(t) v_1 \eta(t) v_2 \eta(t) v_3 
dx dt \right|
\lesssim \lam^{0+} \X 0 \half {v_1}  \X 0 {\half } {v_2}  
\X 0 {0} {v_3} .
\end{equation}
(Strictly speaking, each $\eta$ that appears in \eqref{CaseThreeReduction}
should be replaced by $\eta^{\frac{1}{3}}$ but we abuse the notation with 
the understanding that all smooth cutoff functions are essentially the same
within this analysis.)
We estimate the left-side of \eqref{CaseThreeReduction} using H\"older by
\begin{equation*}
  \lam^\half {{\| \eta(t) v_1 \|}_{L^4_{xt}}}{{\| \eta(t) v_2 \|}_{L^4_{xt}}}
{{\| \eta(t) v_3 \|}_{L^2_{xt}}}
\end{equation*}
In the case under consideration, both $\widehat{v_1}$ and $\widehat{v_2}$
may be assumed to be supported outside $\{ |\xi | \lesssim 1 \}$. In 
such a support region, the prefactor $C(N, \lam )$ appearing in \eqref{cofNlam}
is $\lesssim \lam^{-\frac{1}{4}}$. We apply \eqref{atahalf} to estimate the
$L^4$ norms resulting in the desired estimate
\begin{equation*}
  \lam^{0+} \X 0 \half {v_1 }  \X 0 \half {v_2}  \X 0 0 {v_3} .
\end{equation*}

The preceding discussion established that
\begin{equation}
  \label{firstportion}
  {{\left\| \frac{ {{\langle \xi \rangle}^\half} }{ {{\langle \tau - \xi^3 
\rangle }^\half} {\widehat{ u_1 u_2 }} (\xi ,\tau )} \right\|}_{L^2_\xi 
L^2_\tau}} \lesssim \lam^{0+} \X {-\half} {\half} {u_1}   
\X {-\half} {\half} {u_2}.  
\end{equation}

It remains to control the weighted $L^2_{k} L^1_{\tau}$ portion of
the $Z^{-\half}$ norm to complete the proof of \eqref{bilinearestimate}. 
Since $|\langle \xi  \rangle^{-\half} 
{\widehat{\partial_x (uv) }} ( \xi, \tau )|
\thicksim \langle \xi \rangle^\half | {\widehat{uv}} 
( \xi , \tau )|$, it suffices to
prove that 
\begin{equation}
  \label{toprove}
  {{ \left\| \frac{\langle \xi \rangle ^\half }{ \langle \tau - \xi^3 \rangle} {\widehat{u_1 u_2}} (
\xi, \tau ) \right\|}_{L^2_\xi L^1_\tau}}
\lesssim \lam^{0+} 
\X {-\half } {\half} {u_1}  \X {-\half} {\half} {u_2} .
\end{equation}
The left-side of \eqref{toprove} may be rewritten
\begin{equation}
  \label{rewrittenLone}
  {{\left\| \langle \xi \rangle^\half {{\left\| \frac{1}{ \langle \tau
- \xi^3 \rangle} {\widehat{u_1 u_2 }} (\xi , \tau ) \right\|}_{L^1_\tau}}
\right\|}_{L^2_\xi}}.
\end{equation}

The desired estimate may be reexpressed as
\begin{equation}
  \label{reexpressLone}
  {{\left\| \int\limits_{\xi= \xi_1 + \xi_2 }~  \int\limits_{\tau =
\tau_1 + \tau_2 }  \frac{ \langle \xi \rangle^\half \langle \xi_1 
\rangle^\half \langle \xi_2 \rangle^\half }{\langle \tau - \xi^3 \rangle
\langle \tau_1 - \xi_1^3 \rangle^\half \langle \tau_2 - \xi_2^3 \rangle^\half} 
\widehat{u_1} ( 
\xi_1 , \tau_1 )  \widehat{u_2} ( \xi_2, \tau_2 ) \right\|}_{L^2_\xi 
L^1_{\tau}}}
\end{equation}
\begin{equation*}
\lesssim \lam^{0+} \X 0 0 {u_1}  \X 0 0 {u_2 } . 
\end{equation*}
Recall that the $L^2_\xi$ norm and the various $\xi$-integrations are with
respect to the $\lam$-dependent measure $(d\xi )_\lam$.

Since we may assume our functions have mean zero, we have that $\xi \xi_1
\xi_2 \neq 0$ and the identity
\begin{equation}
\label{arithidentity}
  \tau - \xi^3 = (\tau_1 - \xi_1^3) + (\tau_2 - \xi_2^3 ) - 3 \xi \xi_1 \xi_2
\end{equation}
implies that
\begin{equation}
  \label{denomsstrong}
  \max ( \langle \tau - \xi^3 \rangle, \langle \tau_1 - \xi_1^3 \rangle,  
\langle \tau_1 - \xi_1^3 \rangle ) \gtrsim |\xi \xi_1 \xi_2 |.
\end{equation}
In case $\langle \tau_1 - \xi_1^3 \rangle$ is that maximum, we reduce
matters to \eqref{firstportion}. Indeed, we rewrite \eqref{reexpressLone}
as
\begin{equation*}
  {{\left\| \langle \tau - \xi^3 \rangle^{-\frac{2}{3}}  \langle
\tau - \xi^3 \rangle^{-\frac{1}{3}} \int\limits_{\xi= \xi_1 + \xi_2 }~
\int\limits_{\tau = \tau_1 + \tau_2 } \widehat{u_1} (\xi_1 , \tau_1 )
\widehat{v_2} (\xi_2 , \tau_2 ) \right\|}_{L^2_\xi L^1_{\tau}}}
\end{equation*}
\begin{equation*}
\lesssim \lam^{0+} \X 0 0 {u_1 } \X 0 {\half} {v_2} .
\end{equation*}
Cauchy-Schwarz in $\tau$ (with the observation that $2 ( -\frac{2}{3} ) < -1$)
reduces this case to proving
\begin{equation*}
  {{\left\| \langle \tau - \xi^3 \rangle^{-\third} \int\limits_{\xi= 
\xi_1 + \xi_2 } ~
\int\limits_{\tau = \tau_1 + \tau_2 } \widehat{u_1} (\xi_1 , \tau_1 )
\widehat{u_2} (\xi_2 , \tau_2 ) \right\|}_{L^2_\xi L^2_{\tau}}}
\lesssim \lam^{0+} \X 0 0 {u_1} \X 0 \half {u_2} .
\end{equation*}
Upon rewriting the left-side using duality we see that an $L^4_{xt} 
L^2_{xt} L^4_{xt}$ H\"older application using \eqref{rescaledBourgain}
finishes off this case. The situation when $\langle \tau_2 - \xi_2^3 \rangle$
is the maximum is symmetric so we are reduced to considering the case when
$\langle \tau - \xi^3 \rangle $ is the maximum in \eqref{denomsstrong}.

In the event that
\begin{equation}
  \label{onelittle}
  \langle \tau_1 - \xi_1^3 \rangle \gtrsim |\xi \xi_1 \xi_2 |^{\frac{1}{100}},
\end{equation}
we get a little help from the 1-denominator in \eqref{reexpressLone}. We
cancel $\langle \tau_1 - \xi_1^3 \rangle^{\frac{1}{6}}$ leaving 
$\langle \tau_1 - \xi_1^3 \rangle^\third $ in the denominator and
${{(\langle \xi \rangle \langle \xi_1 \rangle \langle \xi_2 \rangle )}^{
\half -}}$ in the numerator. After the natural cancellation using 
\eqref{denomsstrong}, we collapse to needing to prove
\begin{equation*}
  {{\left\| \langle \tau - \xi^3 \rangle^{-\half - }
\int\limits_{\xi= \xi_1 + \xi_2 } ~
\int\limits_{\tau = \tau_1 + \tau_2 } \widehat{u_1} (\xi_1 , \tau_1 )
\widehat{u_2} (\xi_2 , \tau_2 ) \right\|}_{L^2_\xi L^1_{\tau}}}
\lesssim \lam^{0+} \X 0 {\third} {u_1 } \X 0 {\half} {u_2} .
\end{equation*}
We apply Cauchy-Schwarz in $\tau$ to obtain the upper bound
\begin{equation*}
  {{\left\| 
\int\limits_{\xi= \xi_1 + \xi_2 } ~
\int\limits_{\tau = \tau_1 + \tau_2 } \widehat{u_1} (\xi_1 , \tau_1 )
\widehat{u_2} (\xi_2 , \tau_2 ) \right\|}_{L^2_\xi L^2_{\tau}}}
\end{equation*}
which is controlled as desired using \eqref{rescaledBourgain}.
The case when $\langle \tau_2 - \xi_2^3 \rangle \gtrsim |\xi \xi_1 
\xi_2 |^{\frac{1}{100}} $ is symmetric. 

All that remains is the situation when
\begin{equation}
\label{oncurvebasically}
\langle \tau_i -\xi_i^3 \rangle \ll |\xi \xi_1 \xi_2 |^{\frac{1}{100}}, ~
i=1,2.  
\end{equation}
Recalling \eqref{arithidentity}, we see here that
\begin{equation*}
  \tau - \xi^3 = - 3 \xi \xi_1 \xi_2 + O ( \langle \xi \xi_1 \xi_2 \rangle^{
\frac{1}{100}} ),
\end{equation*}
which we use to restrict $\tau$. After performing the natural cancellation
using \eqref{denomsstrong} on \eqref{reexpressLone}, we wish to show that
\begin{equation}
  \label{omegarestricted}
  {{\left\| \langle \tau - \xi^3 \rangle^{-\half} 
\int\limits_{\xi= \xi_1 + \xi_2 } ~
\int\limits_{\tau = \tau_1 + \tau_2 } \widehat{u_1} (\xi_1 , \tau_1 )
\widehat{u_2} (\xi_2 , \tau_2 ) \chi_{\Omega{(\xi)}} ( \tau - \xi^3 )
\right\|}_{L^2_\xi L^1_\tau}} \lesssim \lam^{0+} \X 0 \half {u_1}
\X 0 \half {u_2} ,
\end{equation}
where the set
\begin{equation*}
\Omega ( \xi ) = \{ \eta \in \R : \eta = - 3 \xi \xi_1 \xi_2 + O 
( \langle \xi \xi_1 \xi_2 \rangle^{\frac{1}{100}} ) ~{\mbox{for any}}~
\xi_1 , \xi_2 \in \Z / \lam ~{\mbox{with}}~ \xi = \xi_1 + \xi_2 \} .
\end{equation*}
We apply Cauchy-Schwarz in $\tau $ to bound the left-side of 
\eqref{omegarestricted} by
\begin{equation*}
  {{\left\| {{\left( \int \langle \tau - \xi^3 \rangle^{-1} \chi_{\Omega(\xi)}
( \tau - \xi^3 ) d\tau \right)}^\half} {{\left\| 
\int\limits_{\xi= \xi_1 + \xi_2 } ~
\int\limits_{\tau = \tau_1 + \tau_2 } \widehat{u_1} (\xi_1 , \tau_1 )
\widehat{u_2} (\xi_2 , \tau_2 ) \right\|}_{L^2_{\tau}}} \right\|}_{L^2_{\xi}}}.
\end{equation*}
The point here is that the characteristic function appearing in the
$\tau$-integrand above sufficiently restricts the region of integration
to prove
\begin{equation}
  \label{laststep}
  {{\left( \int \langle \tau - \xi^3 \rangle^{-1} \chi_{\Omega(\xi)} 
(\tau - \xi^3 ) d\tau \right)}^\half} \lesssim C+ \lam^{0+}
\end{equation}
uniformly in the parameter $\xi$. Note that familiar arguments complete
the proof of \eqref{omegarestricted} (and, hence, \eqref{toprove}) provided
we show \eqref{laststep}.

\begin{remark}
  The condition \eqref{oncurvebasically} restricts the functions 
$\widehat{u_i}$ essentially to the dispersive curve $\{(\xi,\xi^3): \xi
\in \Z / \lam \}$. Suppose for the moment that $\lam = 1$ and we restrict
our attention to those $\xi$ satisfying $|\xi| \thicksim N$. Observe that
the projection of the point set $S_N = \{ (\xi, \xi^3 ) \in \Z_\xi \times 
\R_{\tau}: \xi \in \Z, |\xi|
\thicksim N \}$ onto the $\tau$-axis is a set of $N$ points which are 
$N^2$-separated. Therefore, if we ``vertically thicken'' these points
$O(N^{\alpha})$ for $\alpha \ll 2$, the projected set remains rather 
sparse on the $\tau$-axis. The intuition underlying the proof of 
\eqref{laststep} is that a vertical thickening of 
the set $S_{N_1} + {S_{N_2}}$ also projects onto
a thin set on the $\tau$ axis.
\end{remark}

\begin{lemma}
Fix $\xi \in \Z \backslash \{ 0 \}$.  
For $\xi_1 , \xi_2 \in \Z \backslash \{0\},$
we have for all dyadic $M \geq 1$ that
\begin{equation}
  \label{musize}
  |\{ \mu \in \R : |\mu | \thicksim M, \mu = -3 \xi \xi_1 \xi_2 + O ( \langle 
\xi \xi_1 \xi_2 \rangle^{\frac{1}{100}} ) \} \lesssim M^{1-\delta} 
\end{equation}
for some $\delta > 0$.
\end{lemma}

\begin{proof}
The hypotheses are symmetric in $\xi_1 , \xi_2 $ so we may assume
$|\xi_1 | \geq |\xi_2 |. $ We first consider the situation when $|\xi |
\geq |\xi_1 |$. The expression 
\begin{equation}
  \label{muexpression}
  \mu = -3 \xi \xi_1 \xi_2 + O ( \langle 
\xi \xi_1 \xi_2 \rangle^{\frac{1}{100}} )
\end{equation}
allows us to conclude that 
$|\xi | \lesssim |\mu | \lesssim |\xi|^3$ 
since $\xi_1, \xi_2  \in \Z \backslash 
\{0 \}$ and $|\xi \xi_1 \xi_2 | \lesssim 
|\xi|^3 $. Suppose $|\mu | \thicksim M$ (dyadic) and $|\xi | \thicksim N$ 
(dyadic). We have, for some $p \in [1,3]$, that $M \thicksim N^p$. For $\mu$
to satisfy \eqref{muexpression}, 
$|\xi_1 \xi_2 | \thicksim M^{1 - \frac{1}{p}}$.
We make the crude observation that there are at most $M^{1-\frac{1}{p}}$ 
multiples of $M^{\frac{1}{p}}$ in the dyadic block $\{ |\mu | \thicksim M\}$.
Hence, the set of possible $\mu$ satisfying \eqref{muexpression} must lie
inside a union of $M^{1-\frac{1}{p}}$ intervals of size $M^{\frac{1}{100}}$,
each of which contains an integer multiple of $\xi$. We have then that
\begin{equation*}
  |\{ \mu \in \R : |\mu | \thicksim M, \mu = - 3 \xi \xi_1 \xi_2 + O
( \langle \xi \xi_1 \xi_2 \rangle^{\frac{1}{100}} ) \} |
< M^{{1 - \frac{1}{p} }} M^{\frac{1}{100}} \lesssim M^{\frac{3}{4}},
\end{equation*}
since $1 \leq p \leq 3$.

In case $|\xi | \leq |\xi_1 |$, we must have $|\xi_1 | \lesssim |\mu |
\lesssim |\xi_1 |^3 $ so, if $|\xi_1 | \thicksim N_1$ (dyadic), we 
must have $M \thicksim N_1^p$ for some $p \in [1,3]$ and can repeat the
argument presented above. 
\end{proof}  %of lemma

\begin{remark}
  If we change the setting of the lemma to the case where $\xi, \xi_1, 
\xi_2 \in \Z / \lam \backslash \{0 \}$, 
we have to adjust the conclusion to read
\begin{equation}
\label{lamlemma}
  |\{ \mu \in \R : |\mu | \thicksim M, \mu = - 3 \xi \xi_1 \xi_2 + 
O ( \langle \xi \xi_1 \xi_2 \rangle^{\frac{1}{100}} ) \}| \lesssim
\lam^1 M^{{1-\delta}}, ~\delta>0.
\end{equation}
\end{remark}

We use the lemma to prove \eqref{laststep}. A change of variables leads us
to consider
\begin{equation*}
  \int \langle \mu \rangle^{-1} \chi_{\Omega(\xi ) } ( \mu ) d\mu .
\end{equation*}
We decompose the integration and use \eqref{lamlemma},
\begin{equation*}
  = \int\limits_{|\mu | < \lam^{1000} } \langle \mu \rangle^{-1}
\chi_{\Omega(\xi)} (\mu ) d\mu + \sum\limits_{M: \lam^{1000} < M (dyadic) } ~
\int\limits_{|\mu| \thicksim M} \langle \mu \rangle^{-1} \chi_{\Omega (\xi)}
(\mu ) d\mu
\end{equation*}
\begin{equation*}
  \leq 1000 \log \lam + \sum_{M: \lam^{1000} < M (dyadic)}  
M^{-1} M^{1-\delta} \lam^1 .
\end{equation*}
Finally, we crush $\lam^1$ using the extra decay in $M$ to obtain
\begin{equation*}
  \lesssim \log \lam + \sum_{M (dyadic)} M^{-\frac{\delta}{2}}
\end{equation*}
which proves \eqref{laststep}.
\end{proof}
\subsection{Contraction}

Consider the $\lam$-periodic initial value problem \eqref{lamkdv} with
periodic initial data $\phi \in H^s (0, \lam), ~s \geq - \half.$ We show
first that, for arbitrary $\lam$, this problem is well-posed on a time 
interval of size $\thicksim 1$ provided 
${{\| \phi \|}_{H^{-\half} (0, \lam )}}$ is sufficiently small. Then
we show by a rescaling argument that \eqref{lamkdv} is locally well-posed
for arbitrary initial data $\phi \in H^s (0, \lam)$.

As mentioned before in Remark \ref{meanzeroremark}, 
we restrict our attention to initial data having zero $x$-mean.

Fix $\phi \in H^s (0,\lam ),~s \geq - \half$ and for $w \in Z^{-\half}$ 
define
\begin{equation*}
  \Phi_{\phi} [w] = \eta(t) S_\lam (t) \phi - \eta(t) \int_0^t
S_\lam (t -t') ( \eta(t') w(t') ) dt'.
\end{equation*}
The bilinear estimate \eqref{bilinearestimate} shows that $u \in Y^{-\half}$
implies $\eta(t) \partial_x (u^2) \in Z^{-\half}$ so the (nonlinear)
operator 
\begin{equation*}
\Gamma (u) = \Phi_\phi ( \eta(t) \half \partial_x (u^2))
\end{equation*} 
is defined on $Y^{-\half}$. Observe that $\Gamma (u) = u$ is equivalent,
at least for $t \in [-1,1]$, to \eqref{lamintegraleq}, which is
equivalent to \eqref{lamkdv}.

\noindent{\bf{Claim:}} $\Gamma: ({\mbox{bounded subsets of}} ~Y^{-\half} ) 
 \longmapsto ({\mbox{bounded subsets of}} ~Y^{-\half} ).$ \newline
We estimate
\begin{equation*}
{{\| \Gamma (u) \|}_{Y^{-\half}}} \leq {{\| S_\lam (t) \phi \|}_{Y^{-\half}}}
+ {{\left\| \eta(t) \int_0^t S_\lam ( t- t') ( \eta(t') \half \partial_x 
u^2 (t')) dt' \right\|}_{Y^{-\half}}}.
\end{equation*}
By \eqref{lamhomogestimate} and \eqref{laminhomogestimate}, followed by
the bilinear estimate \eqref{bilinearestimate},
\begin{equation*}
\leq C_1 {{\| \phi \|}_{H^s ( 0, \lam )}} + C_2 {{\| \eta(t) \partial_x
u^2 \|}_{Z^{-\half}}} \leq  C_1 {{\| \phi \|}_{H^s ( 0, \lam )}} + C_2
C_3 \lam^{0+} {{\| u \|}_{Y^{-\half}}^2}
\end{equation*}
and the claim is proven.

Consider the ball
\begin{equation*}
  B = \left\{ 
u \in Y^{-\half} : {{\| w \|}_{Y^{-\half}}} \leq C_4
 {{\| \phi \|}_{H^{-\half}(0, \lam )}} \right\}
\end{equation*}

\noindent
{\bf{Claim:}} $\Gamma $ is a contraction on $B$ if 
${{\| \phi \|}_{H^{-\half}(0,\lam)}}$ is sufficiently small. \newline

We wish to prove that for some $\theta \in (0,1)$, 
\begin{equation*}
  {{\| \Gamma (u) - \Gamma (v) \|}_{Y^{-\half}}} \leq \theta 
{{\| u - v\|}_{Y^{-\half}}}
\end{equation*}
for all $u,v \in B$. Since $u^2 - v^2 = (u+v)(u-v)$, it is not difficult
to see that
\begin{equation*}
  {{\| \Gamma (u) - \Gamma (v) \|}_{Y^{-\half}}} \leq C_2 C_3 \lam^{0+}
( {{\| u \|}_{Y^{-\half}}} + {{\| v \|}_{Y^{-\half}}} ) {{\|  u-v 
\|}_{Y^{-\half}}}.
\end{equation*}
Since $u,v \in B$,
\begin{equation*}
  {{\| \Gamma (u) - \Gamma (v) \|}_{Y^{-\half}}} \leq \lam^{0+}
{{\| \phi \|}_{H^{-\half} (0, \lam )}} {{\|  u-v 
\|}_{Y^{-\half}}}.
\end{equation*}
Hence, for fixed $\lam$, if we take $\phi$ so small that
\begin{equation}
  \label{smalldata}
  \lam^{0+} {{\| \phi \|}_{H^{-\half} (0, \lam )}} \ll 1,
\end{equation}
the contraction estimate is verified.

The preceding discussion establishes well-posedness of \eqref{lamkdv}
on a $O(1)$-sized time interval for any initial data satisfying 
\eqref{smalldata} .

Finally, consider \eqref{lamkdv} with $\lam = \lam_0$ fixed and $\phi
\in H^s ( 0, \lam_0 ), ~s \geq - \half$. This problem is well-posed
on a small time interval $[0, \delta ]$ if and only if the $\sigma$-rescaled
problem
\begin{equation}
  \label{sigmarescaled}
  \left\{
   \begin{matrix}
    \partial_t u_\sigma + \partial_x^3 u_\sigma + \half \partial_x 
u_\sigma^2 =0,& x \in [0, \sigma \lam_0 ] 
         \\
     u_\sigma (x, 0) = \sigma^{-2} \phi (\frac{x}{\sigma} )
   \end{matrix}
\right.
\end{equation}
is well-posed on $[0, \sigma^3 \delta ]$. A simple calculation shows that
\begin{equation*}
  {{\| \phi_\sigma \|}_{H^{-\half}( 0 , \sigma \lam_0 )}}
= \sigma^{-1} {{\| \phi \| }_{H^{-\half } (0, \lam_0 )}}.
\end{equation*}
Observe that
\begin{equation*}
  (\sigma \lam_0 )^{0+} {{\| \phi_\sigma \|}_{H^{-\half} (0 , \sigma \lam_0 )}}
\leq  (\sigma \lam_0 )^{0+}  \sigma^{-1} 
{{\| \phi \| }_{H^{-\half } (0, \lam_0 )}} \ll 1,
\end{equation*}
provided $\sigma = \sigma 
( \lam_0 , {{\| \phi \|}_{H^{-\half} (0 , \lam_0 )}})$ is taken to be
sufficiently large. This verifies \eqref{smalldata} for the problem
\eqref{sigmarescaled} proving well-posedness of \eqref{sigmarescaled}
on the time interval, say $[0, 1]$. Hence, \eqref{lamkdv} is locally
well-posed for $t \in [0, \sigma^{-3} ]$.

\hskip 2in

The preceding discussion reproves the local well-posedness result
for periodic KdV in \cite{KPVBilin}. We record the following simple
variant which will be used in proving the global result for \eqref{TKdV}.
See section 11 of \cite{CKSTTGKdV} for a general interpolation lemma
related to this proposition.

\begin{proposition}
  If $s \geq - \half$, the initial value problem \eqref{TKdV} is locally
well-posed for data $\phi$ satisfying $I \phi \in L^2 ( \T )$. Moreover,
the solution exists on a time interval $[0, \delta]$ with the lifetime
\begin{equation*}
  \delta \thicksim {{\| I \phi \|}_{L^2}^{-\alpha}},
\end{equation*}
and the solution satisfies the estimate
\begin{equation*}
  {{\| I u \|}_{Y^0}} \lesssim {{\| I \phi \|}_{L^2}}.
\end{equation*}
\end{proposition}

\section{Almost conservation and global wellposedness of KdV on $\T$}

This section proves that the 1-periodic initial value problem \eqref{lamkdv}
for KdV is globally well-posed for initial data $\phi \in H^s (\T )$
provided $s \geq -\half$. In particular, we prove Theorem 2. 
The proof is an adaptation of the argument
presented for the real line to the periodic setting.

\subsection{Quintilinear estimate}

The following quinitilinear space-time estimate controls the increment
of the modified energy $E^4_I$ during the lifetime of the local well-posedness
result.

\begin{lemma} 
  \label{Tincrement}
Let $w_i = w_i ( x,t) $ be $\lam$-periodic function 
in $x$ also depending upon $t$. Assume that $\int_0^\lam w_i (x,t) dx = 0$ 
for all $t$. Then
\begin{equation}
  \label{Tquint}
  \left| \int_0^\delta \int_0^\lam \prod_{i=1}^5 w_i (x,t) dx dt \right|
\lesssim \lam^{0+} \prod_{j=1}^3 {{\| w_j \|}_{Y^\half}}  
\X {-\half} {\half} {w_4 }  \X {-\half } {\half } {w_5} .
\end{equation}
\end{lemma}

{\it{Proof (apart from the endpoint):}} 
We group $w_1, w_2, w_3$ together and apply
\eqref{symreduced} to control the left-side by
\begin{equation*}
  \lam^{0+} \X {-\half } {\half} { w_1 w_2 w_3 }  \X {-\half} {\half} {w_4}
  \X {-\half} {\half} {w_5} .
 \end{equation*}
The quintilinear estimate \eqref{Tquint} is thus reduced to proving the
trilinear estimate
\begin{equation}
\label{trilinear}
  \X {-\half } {\half} {w_1 w_2 w_3 } \lesssim \lam^{0+} 
{{\| w_1 \|}_{Y^\half}} {{\| w_2 \|}_{Y^\half}} {{\| w_3 \|}_{Y^\half}}.
\end{equation}

The estimate \eqref{trilinear} is implied by the more general 
fact: For any $s \geq \half$,
\begin{equation}
\label{multilinear}
\X {s-1} \half { \prod_{i=1}^k u_i }
\lesssim \prod_{i=1}^k {{\| u_i \|}_{Y^s}}.
\end{equation}
The multilinear estimate \eqref{multilinear} is proved in 
the forthcoming paper \cite{CKSTTGKdV}. Here we 
indicate the proof for the $k=3$
case of \eqref{multilinear}, namely \eqref{trilinear}, when $s \in (\half, 1]$.
The proof for $s=\half$ in \cite{CKSTTGKdV} 
supplements the discussion presented below with some elementary number theory.
The reader willing to accept \eqref{trilinear} may proceed to Lemma 
\ref{Tapplication}.

The Fourier transform of $\prod_{i=1}^3 {\widehat{u_i}}(x,t)$ equals
\begin{equation}
  \label{trilintransf}
  \int_* \prod_{i=1}^3 \widehat{u_i} ( \xi_i , \tau_i )
\end{equation}
where $\int_*$ denotes an integration over the set where $\xi = \xi_1 +
\xi_2 + \xi_3 ,~ \tau = \tau_1 + \tau_2 + \tau_3 $. We make a case-by-case
analysis by decomposing the left-side of \eqref{trilinear}) into various
regions. We may assume that $\widehat{u_i}~,i = 1,2,3$ are non-negative
$\R$-valued functions.

{\bf{Case 1.}} $\langle \tau - \xi^3 \rangle \lesssim \langle \tau_1 
-\xi_1^3 \rangle.$ \newline
In this case, it suffices to show that
\begin{equation}
  \label{caseonedone}
  \X {s-1} 0 {\prod_{i=1}^3 u_i }  \lesssim  \X s 0 {u_1}  {{\| u_2 \|}_{Y^s}}
{{\| u_3 \|}_{Y^s}} .
\end{equation}
We observe using Sobolev that
\begin{equation*}
  \X {s-1} 0 {\prod_{i=1}^3 u_i }  \lesssim {{\|  {\prod_{i=1}^3 u_i } \|}_{
L^2_t L^{1+}_x} },
\end{equation*}
and then, by H\"older, 
\begin{equation*}
  \lesssim {{\| u_1 \|}_{L^2_t L^{3+}_x}} {{\| u_1 \|}_{L^\infty_t L^{3+}_x}} 
{{\| u_1 \|}_{L^\infty_t L^{3+}_x}} .
\end{equation*}
Finally, using Sobolev again and the embedding $Y^s \subset L^\infty_t H^s_x$,
we conclude that \eqref{caseonedone} holds. Since $\langle \tau_j - \xi_j^3
\rangle \gtrsim \langle \tau - \xi^3 \rangle$ for $j = 2,3$ is symmetric
with the Case 1 defining condition, we may assume that we are in Case 2.

{\bf{Case 2.}} $\langle \tau_i - \xi_i^3 \rangle \ll \langle \tau 
- \xi^3 \rangle $ for $i = 1,2,3.$ \newline
The convolution constraints $\xi = \xi_1 + \xi_2 + \xi_3,~\tau = 
\tau_1 + \tau_2 + \tau_3 $ in this case imply that
\begin{equation*}
  1 \ll \langle \tau - \xi^3 \rangle \thicksim | \xi^3 - (\xi_1^3 + \xi_2^3
+ \xi_3^3 ) |.
\end{equation*}
Therefore, it suffices to show that
\begin{equation}
  \label{casetwodone}
  {{\| \int_* \langle \xi \rangle^{s-1} | \xi^3 - (\xi_1^3 + \xi_2^3
+ \xi_3^3)|^\half 
\widehat{u_1} (\xi_1 , \tau_1 )\widehat{u_2} (\xi_2 , \tau_2 )
\widehat{u_3} (\xi_3 , \tau_3 ) \|}_{L^2_\tau L^2_\xi }}
\end{equation}
\begin{equation*}
  \lesssim \prod_{i=1}^3 \X s \half {u_i} .
\end{equation*}
This estimate may be recast by wiggling the weights and using duality as
\begin{equation}
  \label{dualityrecast}
  \int_* \widehat{f_4} (\xi , \tau ) \frac{ |\xi^3 -  (\xi_1^3 + \xi_2^3
+ \xi_3^3)|^\half}{ \langle \xi \rangle^{1-s} \prod_{i=1}^3 \langle \xi_i 
\rangle^s \langle \tau_i - \xi_i^3 \rangle^\half} 
\widehat{f_1} (\xi_1, \tau_1 )
\widehat{f_2} (\xi_2, \tau_2 )\widehat{f_3} (\xi_3, \tau_3 ) 
\end{equation}
\begin{equation*}
\lesssim \prod_{i=1}^4 {{\| f_i \|}_{L^2_\tau L^2_\xi }}.
\end{equation*}

{\bf{Case IIA.}} $|\xi| \gtrsim |\xi_i |$ for $i = 1,2,3$. \newline

Symmetry allows us to assume $|\xi_1 | \geq |\xi_2 | \geq |\xi_3 |$
and we must have $|\xi | \thicksim |\xi_1 |$. Since $\xi = \xi_1
+ \xi_2 + \xi_3 $, we also have $|\xi^3 - (\xi_1^3 + \xi_2^3 + \xi_3^3 )|
\lesssim |\xi \xi_1 \xi_2 |$ (See \eqref{factfour}). 
Therefore, in this case, the left-side of 
\eqref{dualityrecast} is bounded by
\begin{equation*}
  \int_* \widehat{f_4} ( \xi , \tau ) \frac{ |\xi \xi_1 \xi_2 |^\half}{
\langle \xi \rangle^{1-s} \prod_{i=1}^3 \langle \xi_i \rangle^s \langle
\tau_i - \xi_i^3 \rangle^\half }
\widehat{f_1} ( \xi_1 , \tau_1 )\widehat{f_2} ( \xi_2 , \tau_2 )
\widehat{f_3} ( \xi_3 , \tau_3 )
\end{equation*}
Then, we may bound the preceding by
\begin{equation*}
  \int_* \widehat{f_4}(\xi, \tau) 
\frac{\widehat{f_2} ( \xi_2 , \tau_2 ) }{
\langle \xi_2 \rangle^{s - \half} {\langle 
\tau_2 - \xi_2^3 \rangle^\half}
} 
\frac{\widehat{f_1} ( \xi_1 , \tau_1 )}{\langle 
\tau_1 - \xi_1^3 \rangle^\half}
\frac{
\widehat{f_3} ( \xi_3 , \tau_3 )}{{\langle
\xi_3 \rangle^s } \langle 
\tau_3 - \xi_3^3 \rangle^\half  }
\end{equation*}
and \eqref{dualityrecast} is equivalent to
\begin{equation}
  \label{nextstep}
  {{\| u_1 u_2 u_3 \|}_{L^2_{x,t}}} \lesssim \X 0 \half {u_1}
\X {s-\half} \half {u_2}  \X {s} \half {u_3} .
\end{equation}
We recall from \cite{B1} that $X_{\delta, \half} \subset L^6_{x,t}$
for any $\delta > 0$. Therefore, we validate \eqref{nextstep} using
a H\"older application in $L^4_{x,t} L^6_{x,t} L^{12}_{x,t}$, with the
required $L^{12}_{x,t}$ estimate given by Sobolev and the $L^4_{x,t}$
estimate from \eqref{rescaledBourgain}.

{\bf{Case 2B.}} $|\xi | \ll |\xi_i |$ for $i = 1,2,3.$ \newline
We bound $|\xi^3 - (\xi_1^3 + \xi_2^3 + \xi_3^3 ) | \lesssim |\xi_1
\xi_2 \xi_3 |$ in this region and control the left-side of 
\eqref{dualityrecast} by
\begin{equation*}
  \int_* \widehat{f_4} ( \xi, \tau) \frac{1}{\langle \xi \rangle^{1-s}
\prod_{i=1}^3 \langle \xi_i \rangle^{s - \half} \langle
\tau_i - \xi_i^3 \rangle^\half } \widehat{f_1} ( \xi_1, \tau_1)
\widehat{f_2} ( \xi_2, \tau_2) \widehat{f_3} ( \xi_3, \tau_3)
\end{equation*}
If $s \in (\half, 1]$, we may ignore $\frac{1}{\langle \xi \rangle^{1-s}}$
and finish things off with an $L^2_{x,t} L^6_{x,t} L^6_{x,t} L^6_{x,t}$
H\"older argument using $X_{\delta, \half} \subset L^6_{x,t}$ for any
$\delta >0$.

This completes the proof of \eqref{trilinear} for $s \in (\half, 1]$.

\begin{remark}
  Bourgain has conjectured \cite{B1} that $X_{0 , \half} \subset L^6_{x,t}$.
If this estimate were known, the previous discussion could be substantially
simplified. Our proof of the $s=\half$ case in \cite{CKSTTGKdV} is
partly motivated by an effort to prove this embedding estimate. 
\end{remark}

\begin{lemma} 
\label{Tapplication}
If $m$ is of the form \eqref{particularm} with $s = - \half$ then
\begin{equation}
  \label{TMfiveEstimate}
  \left| \int_0^\delta \Lam_5 (M_5) dt \right| \lesssim \lam^{0+} 
N^{-\frac{5}{2}+}
{{\| Iu \|}_{Y^0}^5} .
\end{equation}
\end{lemma}

The proof is a simple modification of the proof of Lemma \ref{incrementlemma}
with \eqref{Tquint} playing the role of \eqref{quint}. Note also that
$- \frac{3}{4}+$ and $\frac{1}{4}-$ are systematically replaced by
$-\half$ and $\half$ throughout the argument.

\subsection{Rescaling}
Our task is to construct the solution of the 1-periodic \eqref{lamkdv} 
on an arbitrary fixed time interval $[0,T]$. This is equivalent to showing
the $\lam$-rescaled problem with corresponding solution $u_\lam (x,t)
= \lam^{-2} u( \frac{x}{\lam} , \frac{t}{\lam^3} )$ has a solution which
exists on $[0, \lam^3 T]$. The lifetime of the variant local result is 
controlled by ${{\| I \phi \|}_{L^2}}$ and
\begin{equation*}
  {{\| I \phi_\lam \|}_{L^2 (0, \lam )}} \lesssim \lam^{-\frac{3}{2}-s} N^{-s}
{{\| \phi \|}_{H^s (0, \lam )}}.
\end{equation*}
We choose $\lam$ so that
\begin{equation*}
  {{\| I \phi_\lam \|}_{L^2 (0, \lam )}} = \epsilon_0 \ll 1 \implies
\lam \thicksim N^{\frac{-s}{s + \frac{3}{2} }}.
\end{equation*}
This choice guarantees that the local-in-time result for \eqref{lamkdv}
is valid for a time interval of size 1.

\subsection{Almost conservation and iteration}

The local result and Lemma \ref{Tapplication} imply
\begin{equation*}
 | E^4_I ( 1) - E^4_I (0)| \lesssim \lam^{0+} N^{-\frac{5}{2}+} \epsilon_0^5 .
\end{equation*}
Recall that $\lam = \lam(N)$ so we may ignore $\lam^{0+}$ 
by slightly adjusting
$-\frac{5}{2}+$.
Therefore, since by a (modification of) \eqref{twobyfour} 
$E^4_I (t) \thicksim {{\| I \phi_\lam (t) \|}_{L^2}^2}, $
we have that 
\begin{equation*}
  {{\| I \phi_\lam (1) \|}_{L^2}^2} = \epsilon_0 + C \epsilon_0^5  
N^{-\frac{5}{2}+} + O (\epsilon_0^3 )
\end{equation*}
For small $\epsilon_0 $ and large $N$ we see then that ${{\| I \phi_\lam (1)
\|}_{L^2}}$ is also of size $\epsilon_0$. We may iterate the local
result $M$ times until, say,  $E^4_I ( M )$ first exceeds $2 E^4_I (0)$,
that is until
\begin{equation*}
  M N^{-\frac{5}{2}+} \thicksim \epsilon_0 \implies M \thicksim 
N^{\frac{5}{2}-}.
\end{equation*}
The solution of the $\lam (N)$-periodic 
\eqref{lamkdv} is thus extended to the interval $[0, N^{\frac{5}{2}-} ]$. We
now choose $N = N(T)$ such that $N^{\frac{5}{2}-} > [\lam (N)]^3 T 
\thicksim N^\frac{3}{2} T.$ This completes the proof that \eqref{TKdV} is
globally well-posed in $H^{-\frac{1}{2}} ( \T )$. Comments similar
to those presented in \eqref{polyone}-\eqref{polyeight} apply to the
periodic case showing that for our solution of \eqref{TKdV} we have
\begin{equation}
  \label{Tpoly}
  \Hsup {{-\half}} {u(t)} \lesssim t^{\half+} \Hsup {{-\half}} {\phi } .
\end{equation}

\section{Global well-posedness for modified KdV}

The results obtained for KdV are combined with some properties
of the Miura transform \cite{MiuraTransform} (see also 
the survey \cite{MiuraKdVSurvey}, \cite{MiuraKdVSurveyErrata})  to
prove global well-posedness results for modified KdV (mKdV). This section
contains the proofs of Theorems 3 and 4. 
The initial value problem for $\R$-valued mKdV on the line is
\begin{equation}
  \label{mKdVivp}
  \left\{
  \begin{matrix}
    \partial_t u + \partial_x^3 u \pm 6 u^2 \partial_x u =0,& u: \R \times 
[0,T] \longmapsto \R ,
         \\
     u( 0) = u_0.
   \end{matrix}
\right.
\end{equation}
The choice of sign distinguishes between the focussing ($+$) and 
defocussing ($-$) cases. 
This problem is known \cite{KPVCPAM} to be locally well-posed in $H^s$
for $s \geq \frac{1}{4}$. The regularity requirement $s \geq \frac{1}{4}$
is sharp \cite{KPVCounter}. 
We establish global well-posedness of \eqref{mKdVivp}
in the range $s > \frac{1}{4}$ improving the work of Fonseca, Linares and
Ponce \cite{FLP}. 

\subsection{Defocussing case}

Consider the defocussing case of \eqref{mKdVivp}. The {\it{Miura transform}}
of a solution $u$ is the function $v$ defined by
\begin{equation}
  \label{defocusMiura}
  v = \partial_x u + u^2.
\end{equation}
A calculation shows that $v$ solves
\begin{equation}
  \label{MiuraKdVivp}
  \left\{
   \begin{matrix}
    \partial_t v + \partial_x^3 v - 6 v \partial_x v =0,& v: \R \times 
[0,T] \longmapsto \R ,
         \\
     v( 0) = v_0.
   \end{matrix}
\right.
\end{equation}

Suppose the initial data $u_0$ for \eqref{mKdVivp} is in $H^s,~\frac{1}{4} < s
< 1$. We show that $v_0$ is in $H^{s-1}$. 

\begin{lemma}
\label{tenonelemma}
  If $u_0 \in H^s,~\frac{1}{4} < s < 1$, then $v_0 = \partial_x u_0 + u_0^2
\in H^{s-1}$.
\end{lemma}

\begin{proof}
  \begin{eqnarray*}
    \Lp 2 { D^{s-1} v_0 } &\leq &\Lp 2 {D^{s-1} \partial_x u_0} 
+ \Lp 2 {D^{s-1} ( u_0^2 ) } \\
   & \leq & \Lp 2 {D^s u_0 } + \Lp 2 {u_0^2 } \\
   &  \leq& \Lp 2 {D^s u_0 } + {{\| u_0 \|}_{L^4}^2 }  \\
& \lesssim  &
\Lp 2 {D^s u_0 }
+ {{\| u_0 \|}_{H^s }^2}.
   \end{eqnarray*}
\end{proof}

The lemma verifies that the initial data $v_0$ for the KdV equation 
\eqref{MiuraKdVivp} is in $H^{s-1}$ and $-\frac{3}{4} < {s-1}$ since
$\frac{1}{4} < s$. Therefore, the global well-posedness result for KdV
just established
applies to \eqref{MiuraKdVivp} and we know that the solution $v$ 
exists for all time and satisfies
\begin{equation}
  \label{polybound}
  \Hsup {{s-1}} {v(t)} \lesssim (1+ |t|)^C
\end{equation}
for some constant $C$. We exploit this polynomial-in-time bound for
KdV solutions to control $\Hsup s {u(t)}$ using the Miura transform.

Note that 
\begin{equation*}
  \Hsup s {u(t)} \lesssim \Lp 2 {u(t)} + \Hsup {s-1} {\partial_x u (t) }.
\end{equation*}
Since our mKdV solution satisfies $L^2$-mass conservation, $\Lp 2 {u(t)}
= \Lp 2 {u_0}$, it suffices to control $\Hsup {{s-1}} {\partial_x u(t)}$
to control $\Hsup s {u(t)}$. By \eqref{defocusMiura},
\begin{equation*}
  \Hsup {{s-1}} {\partial_x u(t) } \lesssim \Hsup {{s-1}} {v(t) } +
\Hsup {{s-1}} {u^2 (t) }.
\end{equation*}
Using \eqref{polybound},
\begin{equation*}
  \Hsup {{s-1}} {\partial_x u(t) } \lesssim (1+|t|)^C +
\Hsup {{s-1}} {u^2 (t) }.
\end{equation*}
Summarizing, we have
\begin{equation}
  \label{almostpinched}
  \Hsup s {u(t)} \lesssim 1 + (1+|t|)^C + \Hsup {{s-1}} {u^2 (t) }.
\end{equation}

\begin{lemma}
\label{tentwolemma}
  Assuming that $\Lp 2 {u(t)} \lesssim 1$ and $\frac{1}{4} < s < 1$, there
exists an $\epsilon > 0$ such that
\begin{equation}
  \label{claimtwo}
  {{\| u^2 (t) \|}_{H^{s-1}}} \lesssim {{\| u(t) \|}_{H^s}^{1-\epsilon}}.
\end{equation}
\end{lemma}

Assuming the lemma for a moment, observe that 
combining \eqref{claimtwo} and \eqref{almostpinched} implies a 
polynomial-in-time upper bound on $\Hsup s {u(t)}$ giving global 
well-posedness of defocussing mKdV. We now turn to the proof of 
\eqref{claimtwo}.
\begin{proof}
  We first consider the case when $\half + \frac{1}{1000} 
< s < 1$. The Sobolev 
estimate in one dimension
\begin{equation}
  \label{GagNir}
  \Lp q w  \lesssim \Lp p {D^\sigma w }; ~~~~\frac{1}{q} = \frac{1}{p} - 
\frac{\sigma}{1}
\end{equation}
is applied with $w = D^{s-1} (u^2),~\sigma = (1-s), ~q=2$, yielding for
\begin{equation}
\label{pcondition}
\frac{1}{p} = \frac{3}{2} - s
\end{equation}
that
\begin{equation*}
  {{\| D^{s-1} (u^2) \|}_{L^2}} \lesssim {{\| u^2 \|}_{L^p}}
\end{equation*}
We continue the estimate by writing ${{\| u^2 \|}_{L^p}} = 
{{\| u \|}_{L^{2p}}}$ and using Sobolev to get
\begin{equation}
\label{continue}
  \leq {{\| u \|}_{L^{2p}}^2} \leq {{\| u \|}^2_{H^{\sigma(p)}}}; ~\sigma (p )
= \half - \frac{1}{2p}.
\end{equation}
Finally, we interpolate $H^{\sigma(p)}$ between $H^0 = L^2 $ and $H^s$ to
obtain
\begin{equation*}
  {{\| D^{s-1} (u^2) \|}_{L^2}} \lesssim {{\| u \|}_{L^2}^{2 (1-\theta)}}
{{\| u \|}_{H^s}^{2 \theta }}
\end{equation*}
where $\theta = \frac{1}{s} \sigma ( p )$. Using \eqref{pcondition} and
\eqref{continue}, we can simplify to find $\sigma(p) = \frac{2s -1}{4}$
and $2 \theta = 1- \frac{1}{2s}$. Since $\Lp 2 {u} \lesssim 1$, we observe
that \eqref{claimtwo} holds in case $\half + \frac{1}{1000} < s < 1.$

In case $\frac{1}{4} < s \leq \half + \frac{1}{1000}$, 
we begin with a crude step
by writing 
\begin{equation*}
  {{\| D^{s-1} (u^2) \|}_{L^2}} \lesssim  {{\| D^{s- \frac{2}{3} } 
(u^2) \|}_{L^2}} .
\end{equation*}
Modifying the steps in the previous case, we have
\begin{equation*}
  {{\| D^{s - \frac{2}{3}} (u^2) \|}_{L^2}} \leq {{\| u \|}_{L^{2p}}^2}; ~
\frac{1}{p} = \half + \frac{2}{3} -s .
\end{equation*}
Then, by Sobolev and interpolation,
\begin{equation*}
  {{\| u \|}_{L^{2p}}^2} \lesssim {{\| u \|}_{H^{\sigma{(p)}}}^2} \lesssim
{{\| u \|}_{L^2}^{2(1-\theta)}} {{\| u \|}_{H^s}^{2 \theta}} ,
\end{equation*}
where
\begin{equation*}
  \theta = \frac{1}{s} \sigma(p), ~ \sigma (p ) = \half - \frac{1}{2p} =
\frac{s}{2} - \frac{1}{12}.
\end{equation*}
It is then clear that for $s \in (\frac{1}{4}, \half + \frac{1}{1000}]$, 
we have 
$2 \theta = 1 - \frac{1}{6s} = 1 - \epsilon$ for
an appropriate $\epsilon > 0$ as claimed.

\end{proof}

\subsection{Focussing case}

In the focussing case of $\R$-valued 
modified KdV, the Miura transform has a different
form
\begin{equation}
  \label{focusMiura}
  v = \partial_x u + i u^2 .
\end{equation}
The function $v$ solves the {\it{complex}} KdV initial value problem
\begin{equation}
  \label{CKdVivp}
  \left\{
  \begin{matrix}
    \partial_t v + \partial_x^3 v -i 6 v \partial_x v =0,& v: \R \times 
[0,T] \longmapsto \C ,
         \\
     v( 0) = v_0.
   \end{matrix}
\right.
\end{equation}

Since the solution $u(t)$ of focussing modified KdV is $\R$-valued
and derivatives are more costly than squaring in one dimension, we take the
perspective that $v$ is ``nearly $\R$-valued''. The variant local
result for \eqref{CKdVivp} has an existence interval determined by
${{(\int |I v_0 |^2 dx )}^\half}$. However, \eqref{CKdVivp} does not
conserve ${{(\int |I v (t) |^2 dx )}^\half}$ but instead (almost) conserves
${{|\int (Iv(t))^2 dx |}^\half}$. An iteration argument showing global
well-posedness may proceed if we show that
\begin{equation*}
{{\left|\int (Iv(t))^2 dx \right|}^\half} {~{\mbox{controls}}~}
{{\left(\int |I v (t) |^2 dx \right)}^\half} 
\end{equation*}
for functions $v$ of the form given by the Miura transform \eqref{focusMiura}.

Observe that
\begin{equation}
\label{modulus}
{{\left(\int |Iv |^2 dx \right)}^\half }
= {{\left( \int (I u_x )^2 + (I(u^2))^2 dx \right)}^\half},
\end{equation}
\begin{equation}
\label{square}
{{\left| \int (Iv )^2 dx \right| }^\half } = {{\left| \int (I u_x )^2 -
(I ( u^2))^2 dx + 2i \int (I u_x ) (I (u^2)) dx \right|}^\half}.
\end{equation}

\begin{lemma}
Assuming that ${{\| u(t) \|}_{L^2}} \leq C $ for all $t$, and $\frac{1}{4}
< s < 1$, there exists an $\epsilon > 0$ such that
\begin{equation}
{{\| I (u^2 ) \|}_{L^2}} \leq \epsilon {{ \| I u_x \|}_{L^2}} + C.
\label{squarelowerorderestimate}
\end{equation}
\label{squarelowerorderlemma}
\end{lemma}

If we take the lemma for granted, we deduce from \eqref{modulus} and
\eqref{square} that 
\begin{equation}
\label{modlikesquare}
{{\left| \int |I v(t) |^2 dx \right|}^\half} {\mbox{is bounded}}
\iff {{\left( 
\int (I v(t) )^2 dx \right)}^\half} {\mbox{is bounded}} 
\end{equation}
\begin{equation*}
\iff
{{\left| \int 
(I u_x (t) )^2 dx 
\right|}^\half} {\mbox{is bounded}}.
 \end{equation*}

The equivalence \eqref{modlikesquare} links the quantity determining the
length of the local existence interval to an almost conserved quantity.
Consequently, \eqref{CKdVivp} is GWP and ${{| \int (I v(t))^2 dx|}^\half}$
is polynomially bounded in $t$. Since $\Lp 2 {u(t)} \leq C$ for solutions
of focussing modified KdV, the equivalence \eqref{modlikesquare} implies
${{\| u(t) \|}_{H^s }}, s > \frac{1}{4}$ is polynomially bounded
in $t$. Therefore, focussing mKdV is globally well-posed in $H^s (\R ),~
s> \frac{1}{4},$ {\it{provided}} we prove the lemma above.

\begin{proof}[Proof of Lemma \ref{squarelowerorderlemma}:] 

We use duality and
rewrite the expression to be controlled as
\begin{equation}
\label{tensixteen}
  \int f( \xi_1 + \xi_2 ) m( \xi_1 + \xi_2 ) \widehat{u_1} ( \xi_1 )
\widehat{u_2} ( \xi_2 ) d\xi_1 d\xi_2 
\end{equation}
where $m( \xi ) \thicksim 1$ when $|\xi | \lesssim N$ and $m( \xi )
\thicksim N^{-(s-1)} |\xi |^{s-1} $ when $|\xi | > N$ and we have
relaxed to the bilinear situation. The function $f$ is introduced
to calculate the norm using duality so $\Lp 2 f \leq 1$. We may assume
that $\widehat{u_j}$ is nonnegative.
Symmetry allows us to assume
$|\xi_1 | \geq |\xi_2 |.$

{\bf{Case 1.}} $|\xi_1 | \lesssim N \implies m( \xi_1 + \xi_2 ) \thicksim 1.$
\newline
In this case, $I$ acts like the identity operator and the task is to 
control ${{\| u \|}_{L^4}^2}$. By Sobolev and interpolation,
\begin{equation*}
  {{\| u \|}_{L^4}^2} \leq {{\| u \|}_{L^2}} {{\| u \|}_{{\dot{H}}^\half }}
\leq {{\| u \|}_{L^2}^{\frac{3}{2}}} {{\| u \|}_{{\dot{H}}^1}^{\frac{1}{2}}}
\end{equation*}
and we observe that, in this case, 
\begin{equation*}
  {{\| I( u^2 ) \|}_{L^2}} \leq C + \epsilon 
{{\| I (\partial_x u) \|}_{L^2}}.
\end{equation*}

{\bf{Case 2.}} $|\xi_1 | \gg N.$
\newline
{\bf{Case 2.A.}} $|\xi_1 + \xi_2 | \lesssim N \implies |\xi_1 | \thicksim 
|\xi_2 |,~ m( \xi_1 + \xi_2 ) \thicksim 1.$
\newline
We decompose the factors dyadically by writing
\begin{equation*}
  \int f( \xi_1 + \xi_2 ) \widehat{u_1 } ( \xi_1 ) \widehat{u_2 } ( \xi_2 )
d\xi_1 d \xi_2  \thicksim \sum_{N_1 \thicksim N_2 \gg N}
\int f( \xi_1 + \xi_2 ) {\widehat{u_{N_1}}} ( \xi_1 ) {\widehat{u_{N_2}}} 
( \xi_2) d\xi_1 d\xi_2 .
\end{equation*}
We focus on a particular dyadic interaction term
\begin{equation}
 \label{dyadicinteraction}
  \int_{\{ |\xi_i |\thicksim N_i \}} f( \xi_1 + \xi_2 ){\widehat{u_{N_1}}} 
( \xi_1 ) {\widehat{u_{N_2}}} 
( \xi_2) d\xi_1 d\xi_2 . 
\end{equation}
\begin{equation*}
  \thicksim  \int_{\{ |\xi_i | \thicksim N_i \}} f( \xi_1 + \xi_2 )
m( \xi_1 ) \xi_1 {\widehat{u_{N_1}}} 
( \xi_1 ) \frac{1}{m(\xi_1 ) \xi_1 } {\widehat{u_{N_2}}} 
( \xi_2) d\xi_1 d\xi_2 . 
\end{equation*}
In the Case 2.A. region, $m( \xi_1 ) \xi_1  \thicksim N_1^s N^{-(s-1)}$.
We make this substitution and apply Cauchy-Schwarz in $\xi_2$ to observe
\begin{equation*}
  \lesssim {{\| f \|}_{L^2_{\{ |\xi | \lesssim N \}}}}
{{\| {\widehat{u_{N_2}}} \|}_{L^2}} \frac{1}{N_1^s}  {{\| |\xi_1|^s
{\widehat{u_{N_1}}}(\xi_1)  \|}_{L^1}}.
\end{equation*}
Multiplying through by $1 = \frac{N^{-(s-1)}}{N^{-(s-1)}}$ leads to
\begin{equation*}
  \lesssim \Lp 2 f  \Lp 2 {\widehat{u_{N_2}}} \frac{N^{s-1}}{N_1^{s-\half}}
\Lp 2 {I \partial_x u_{N_1}}
\end{equation*}
\begin{equation*}
  \ll \Lp 2 {I \partial_x u_{N_1}}, ~{\mbox{provided}}~\half< s < 1,
\end{equation*}
since $N_1 > N \gg 1.$ Of course we can sum over the dyadic scales
and retain the claim.

It remains to establish the claim for the Case 2.A. region 
when $\frac{1}{4} < s \leq \frac{1}{2}$.
We rewrite the expression \eqref{dyadicinteraction} differently as
\begin{equation*}
  \int \frac{ f( \xi_1 + \xi_2 ) }{{{(1+ |\xi_1 + \xi_2 |)}^{\half +}}}
{{(1+ |\xi_1 + \xi_2 |)}^{\half+}} {{|\xi_1 |}^s} {\widehat{u_{N_1}}}
( \xi_1 )  \frac{1}{|\xi_1|^s } {\widehat{u_{N_2}}} ( \xi_2 ) d\xi_1
d\xi_2 . 
\end{equation*}
Defining $\widehat{D} (\xi ) =  \frac{ f( \xi_1 + \xi_2 ) }{{{(1+ |\xi_1 + \xi_2 |)}^{\half+}}}$, $\widehat{U_{N_1}}(\xi_1 ) = |\xi_1|^s \widehat{u_{N_1}} 
( \xi_1),$ we observe that the preceding expression is controlled by
\begin{equation*}
  N^{\half+} N_1^{-s} 
\langle \widehat{D} , \widehat{U_{N_1}} * \widehat{u_{N_2}} \rangle.
\end{equation*}
We apply H\"older to estimate by
\begin{equation*}
  N^{\half +} N_1^{-s} \Lp \infty D  \Lp 2 {U_{N_1}} \Lp 2 {u_{N_2}}
\end{equation*}
and Sobolev implies $\Lp {\infty} D \leq \Lp 2 f$. Rewriting this gives
\begin{equation*}
  N^{\half+} \frac{1}{N_1^s} \Lp 2 f  \Lp 2 { |\xi_1 |^s {\widehat{u_{N_1}}}
( \xi_1 ) } \Lp 2 {u_{N_2}}.
\end{equation*}
Multiplying by $\frac{N^{-(s-1)}}{N^{-(s-1)}}$ yields
\begin{equation*}
  \left( \frac{N}{N_1} \right)^s \frac{1}{N^{\half-}} \Lp 2 f
\Lp 2 {u_{N_2}} \Lp 2 {I \partial_x u_{N_1}}
\end{equation*}
and the prefactor vanishes as $N \rightarrow \infty$, proving the claimed
estimate.

{\bf{Case 2.B.}} $|\xi_1 + \xi_2 | \gg N$. 
\newline
We return to \eqref{tensixteen} with a particular dyadic interaction.
We multiply and divide through by $m(\xi_1) \xi_1$ and use the
definition of $m(\cdot )$ to see
\begin{equation*}
  \int |\xi_1 + \xi_2 |^{- (s-1) } N^{-(s-1)} f(\xi_1 + \xi_2 )
~m(\xi_1 ) \xi_1 \widehat{u_{N_1}} (\xi_1 ) ~\frac{1}{N_1^s N^{-(s-1)}}
\widehat{u_{N_2}} (\xi_2 ) ~d\xi_1 d\xi_2 
\end{equation*}
and, since $|\xi_1 + \xi_2 | \lesssim N_1$,
\begin{equation*}
  \lesssim \int f(\xi_1 + \xi_2 ) ~m(\xi_1 ) \xi_1  
\widehat{u_{N_1}} (\xi_1 ) ~\frac{1}{N_1} \widehat{u_{N_2}} (\xi_2 ) 
~d\xi_1 d\xi_2 .
\end{equation*}
We estimate $L^2~L^2~L^1$ and go back to $L^2$ in the last factor. The
return to $L^2$ introduces $N_2^\half$ which is compensated by $N_1^{\half}$
in the denominator and we retain enough decay in $N_1$ to sum.
\end{proof}

\subsection{Modified KdV on $\T$ }

The Lemmas \ref{tentwolemma} and \ref{squarelowerorderlemma} naturally extend
to the $\lam$-periodic setting. These results link the polynomial-in-time
upper bound \eqref{Tpoly} for solutions of the $\lam$-periodic initial value
problem \eqref{lamkdv} for KdV to a polynomial-in-time upper bound on
${{\| u(t) \|}_{H^\half (\T )}}$ for solutions of the $\lam$-periodic initial
value problem for mKdV, implying Theorem 4.

% \bibliographystyle{plain}
% \bibliography{\homedir/TeX/master}
% %\bibliography{/home/colliand/Math/Papers/Biblio/master}

\enddocument